\documentclass[11pt]{article}

\usepackage{latexsym}
\usepackage{amssymb}
\usepackage{amsthm}
\usepackage{amscd}
\usepackage{amsmath}
\usepackage{tikz}
\usetikzlibrary{patterns}
\usepackage{mathabx}
\usepackage{stmaryrd}
\usepackage{shuffle}
\usepackage{comment}
\usepackage{dsfont}
\usepackage[margin=1in]{geometry}

\usepgflibrary{arrows}

\newtheorem{theorem}{Theorem}[section]

\newtheorem{lemma}[theorem]{Lemma}

\newtheorem{proposition}[theorem]{Proposition}
\newtheorem{corollary}[theorem]{Corollary}

\theoremstyle{definition}

\newtheorem*{runningexample}{Running Example}
\newtheorem*{examples}{Examples}
\newtheorem{remark}[theorem]{Remark}

\numberwithin{equation}{section}

\newcommand{\CC}{\mathbb{C}} 
\newcommand{\FF}{\mathbb{F}}

\newcommand{\ZZ}{\mathbb{Z}}

\newcommand{\cB}{\mathcal{B}}

\newcommand{\cH}{\mathcal{H}}

\newcommand{\fkut}{\mathfrak{ut}}

\newcommand{\Ind}{\mathrm{Ind}}

\newcommand{\Res}{\mathrm{Res}} 
\newcommand{\Inf}{\mathrm{Inf}}

\newcommand{\UT}{\mathrm{UT}}

\newcommand{\Def}{\mathrm{Def}}

\newcommand{\Sym}{\mathrm{Sym}}

\newcommand{\scf}{\mathrm{scf}}

\newcommand{\cf}{\mathrm{cf}}

\newcommand{\dd}{\displaystyle}
\newcommand{\scs}{\scriptstyle}
\newcommand{\scscs}{\scriptscriptstyle}

\newcommand{\spanning}{\textnormal{-span}}

\def\adots{\mathinner{\mkern2mu\raise0pt\hbox{.}  
\mkern2mu\raise4pt\hbox{.}\mkern1mu
\raise7pt\vbox{\kern7pt\hbox{.}}\mkern1mu}}

\newcommand{\Cl}{\mathtt{Cl}}
\newcommand{\Ch}{\mathtt{Ch}}
\newcommand{\reg}{\mathbf{reg}}

\newcommand{\NCSym}{\mathrm{NCSym}}
\newcommand{\NSym}{\mathrm{NSym}}

\newcommand{\battery}[1]{\curlyvee_{\hspace*{-.5ex}A}}

\newcommand{\lc}{\mathrm{lc}}
\newcommand{\bc}{\mathrm{bc}}

\newcommand{\llc}{\mathrm{llc}}

\newcommand{\Down}{\mathrm{Dn}}
\newcommand{\TF}{\mathrm{TF}}

\newcommand{\QSym}{\mathrm{QSym}}
\newcommand{\sign}{\mathrm{sgn}}
\newcommand{\toggle}{\mathrm{tgl}}
\newcommand{\prim}{\mathrm{prim}}
\newcommand{\binaryzero}{\mathrm{bn}_0}
\newcommand{\binaryone}{\mathrm{bn}_1}
\newcommand{\Splt}{\mathrm{Splt}}
\newcommand{\Fusd}{\mathrm{Fusd}}
\newcommand{\twst}{\mathrm{twst}}
\newcommand{\desc}{\mathrm{desc}}
\newcommand{\MR}{\mathrm{MR}}
\newcommand{\vdotdotv}{\hspace{-.5ex}\shortmid\hspace{-.9ex}\cdot \cdot\hspace{-.9ex}\shortmid\hspace{-.5ex}}
\newcommand{\dotvdot}{\cdot\hspace{-.75ex} \shortmid\hspace{-.75ex}  \cdot}
\newcommand{\sdotvdot}{\cdot \shortmid  \cdot}
\newcommand{\svdotdotv}{\shortmid\cdot \cdot\shortmid}
\newcommand{\all}{\bullet}
\newcommand{\dual}{*}

\makeatletter
\renewcommand{\@makefnmark}{\mbox{\textsuperscript{}}}
\makeatother

\allowdisplaybreaks[1]

\title{Hopf structures in the representation\\ theory of direct products}

\date{}
\author{Farid Aliniaeifard and Nathaniel Thiem}

\begin{document}

\maketitle

\begin{abstract}
Combinatorial Hopf algebras give a linear algebraic structure to infinite families of combinatorial objects, a technique further enriched by the categorification of these structure via the representation theory of families of algebras.  This paper examines a fundamental construction in group theory, the direct product, and how it can be used to build representation theoretic Hopf algebras out of towers of groups.  A key special case gives us the noncommutative symmetric functions $\NSym$, but there are many things that we can say for the general Hopf algebras, including the structure of their character groups and a formula for the antipode.{\let\thefootnote\relax}
\end{abstract}

\section{Introduction}

Combinatorial Hopf algebras have been an essential tool used to study infinite families of combinatorial objects even before they were formally defined in \cite{ABS}.   However, in our favorite examples, the true potential becomes far more apparent when we have a representation theoretic realization of the Hopf algebra, such as the fundamental connection between the Hopf algebra of symmetric functions $\Sym$ and the representation theory of the symmetric groups \cite{Mac}.    There are numerous approaches to such ``categorifications" of combinatorial Hopf algebras in the literature, whether by $0$-Hecke algebras \cite{KT}, $0$-Hecke--Clifford algebras \cite{BHT} or quantum groups \cite{KL}.  This paper takes the approach from \cite{Aim,ATca,ATmr}, applying recent ideas in the representation theory of groups to categorify combinatorial of Hopf algebras via towers of groups.

Supercharacter theories of finite groups are fundamental to this approach.  The idea was originally introduced by Andr\`e (e.g.~\cite{An}) to study the representation theory of wild groups, such as the finite groups of unipotent uppertriangular matrices. Diaconis and Isaacs later formalized the concept in \cite{DI}.   While understanding all the supercharacter theories of groups is generally hard (even for abelian groups), one can interpret a choice of supercharacter theory as a way to specify which aspects of the representation theory one would like to focus on. Each supercharacter theory of a group $G$ replaces the space of class functions $\cf(G)$ with a subspace of \emph{superclass functions} $\scf(G)$.  For example, a key case in this paper gives the two-dimensional subspace
\begin{equation}\label{TrivialTheory}
\scf(G)=\{\psi:G\rightarrow\CC\mid \psi(g)=\psi(h), \text{ if $g\neq 1$ and $h\neq 1$}\}.
\end{equation}
Effectively, this theory only distinguishes between what is trivial and not trivial.

The paper \cite{Aim} first demonstrated that one could use the supercharacter theory of a tower of groups to construct a representation theoretic realization of a Hopf algebra; in that case, the authors categorified the Hopf algebra of symmetric functions in non-commuting variables $\NCSym$ via the representation theory of the finite unipotent upper-triangular matrices.  While the isomorphism given in the paper is natural, it lacks the canonical nature of the symmetric group ($\Sym$) case.  One underlying philosophy of this paper is that this ``deficiency" may turn out to be a feature of the theory.

Here we consider the fundamental group theory construction of taking direct products of groups, 
\[G^{n-1}=\underbrace{G\times G\times \cdots \times G}_{\text{$n-1$ factors}}\times \{1\}\]
and builds a tower of groups.  For the representation theory, we fix a single supercharacter theory for $G$ and extend it to the direct product groups to get
\[
\scf(G^{n-1})\cong \scf(G)\otimes \cdots \otimes \scf(G)\otimes \scf(\{1\}).
\]
To obtain a Hopf structure, we develop product and coproduct functors that rely on the inclusion/quotent given by
\[
G^{m-1}\times G^{n-1} \begin{array}{c} \overset{\text{inclusion}}{\longrightarrow} \\[-.25cm] \underset{\text{projection}}{\longleftarrow} \end{array} G^{m+n-1},
\]
allowing us to use both induction and inflation.  The resulting family of Hopf algebras
\[
\{\cH_{(\iota, \alpha,\beta)} \mid \iota,\alpha,\beta\in \scf(G), \langle\iota,\alpha\rangle=\langle \iota, \beta\rangle=1\}
\]
 depend on supercharacter theory of $G$, and all have a structure motivated by the Hopf algebra of noncommutative symmetric functions $\NSym$.  This connection becomes an explicit isomorphism in the special case when we take the theory (\ref{TrivialTheory}) of $G$, but the isomorphisms vary widely. 
 
In fact, we show that for this supercharacter theory,
 \[
 \cH_{(\iota,\alpha,\beta)}\cong\NSym,
 \]
 for any choice of $\iota$, $\alpha$ and $\beta$, but each choice actually makes certain internal constructions more natural.   Some other highlights include
\begin{description}
\item[Section \ref{NSym}.]  An analysis of how to use $\alpha$, $\beta$ and $\iota$ to construct our favorite bases of $\NSym$ inside of $\cH_{(\iota,\alpha,\beta)}$,
\item[Section \ref{Embeddings}.] Realizations of the standard embeddings of $\NSym$ into $\NCSym$ and the Malvenuto--Reutenauer Hopf algebra $\MR$ using representation theoretic functors by varying our choice of $\iota$, $\alpha$ and $\beta$.
\item[Section \ref{CharacterGroupSection}.] An analysis of the structure of the character group of $\cH_{(\iota,\alpha,\beta)}$ as a combinatorial Hopf algebra, including a characterization of when a function in  the dual $\cH_{(\iota,\alpha,\beta)}^*$ is in fact a character, and a formula for the inverse of a character, depending on $\alpha+\beta$.
\item[Section \ref{Antipode}.] A formula for the antipode in $\cH_{(\iota,\alpha,\beta)}$. 
\end{description} 

The paper is organized as follows.  Section \ref{Preliminaries} reviews supercharacter theories, combinatorial Hopf algebras, and some of the basic combinatorics used in this paper.  Section \ref{HopfAlgebraConstruction} defines the Hopf algebras $\cH_{(\iota,\alpha,\beta)}$ using some representation theoretic functors.  Section \ref{SectionExamples} specializes to the case where we recover $\NSym$, specifying some explicit isomorphisms and connecting this example to the Hopf algebras in  \cite{Aim} and \cite{ATmr}.  Section \ref{CombinatorialHopfAlgebras} returns to the general setting and concludes with some analyses of the character groups of $\cH_{(\iota,\alpha,\beta)}$ and the antipode.

\vspace{.5cm}

\noindent\textbf{Acknowledgements.} Thiem was supported by Simons Foundation collaboration grant 426594. 

\section{Preliminaries} \label{Preliminaries}

We begin by a brief review of supercharacter theories, keeping in mind that most of our main results use only unsubtle versions of the theory.  We then set our notation for integer and set compositions and conclude with a review of combinatorial Hopf algebras, along with some favorite examples.

\subsection{Supercharacter theories}

A \emph{supercharacter theory} $(\Ch,\Cl)$ of a finite groups $G$  is a set partition $\Cl$ of $G$ with a basis of orthogonal characters $\Ch$ for the space
\[
\scf(G)=\{\psi:G\rightarrow \CC\mid \text{$\psi$ constant on the blocks of $\Cl$}\}
\]
such that the regular character $\reg\in\scf(G)$ and the trivial character $\mathds{1}\in \Ch$, where orthogonality is determined by the restriction of the usual inner product  $\langle\cdot,\cdot\rangle:\cf(G)\otimes \cf(G)\longrightarrow \CC$ given by
\[
\langle\psi,\gamma\rangle=\frac{1}{|G|}\sum_{g\in G} \psi(g)\overline{\gamma(g)}
\]
to $\scf(G)$.   We generally refer to the elements of $\Ch$ as \emph{supercharacters} and the blocks of $\Cl$ as \emph{superclasses}. 

Any such a theory comes equipped with two canonical bases of $\scf(G)$: the supercharacters $\Ch$ and the set superclass identifier functions
\[
\{\delta_A\mid A\in \Cl\}\quad \text{where}\quad \delta_A(g)=\left\{\begin{array}{ll}
1 & \text{if $g\in A$,}\\
0 & \text{if $g\notin A$.}
\end{array}\right.
\]

 The main example for this paper is the (trivial?) case where
\[
\Cl=\{\{1\},G-\{1\}\}\quad \text{and}\quad \Ch=\{\mathds{1},\reg-\{\mathds{1}\}\} \quad \text{with}\quad \scf(G)=\CC\spanning\{\mathds{1},\reg\},
\]
and much of the paper can be understood with only this example and the following basic generalization to direct products of groups.  

If $G$ is a group and 
\[G^{n-1}=G\times G\times \cdots \times G\times \{1\}\]
is a direct product, then for any supercharacter theory $(\Ch_G,\Cl_G)$ of $G$, we define a supercharacter theory of $G^{n-1}$ by
\begin{equation} \label{DirectProductSupercharacterTheory}
\begin{split}
\Cl &=\{(A_1,\ldots,A_{n-1},\{1\})\mid A_1,\ldots, A_{n-1}\in \Cl_G\}
\\
 \Ch &=\{(\chi_1,\ldots,\chi_{n-1},\chi^{()})\mid \chi_1,\ldots, \chi_{n-1}\in \Ch_G\},
\end{split}
\end{equation}
where $\chi^{()}$ is the trivial character of the trivial group $\{1\}$.  Then the corresponding function space satisfies
\[
\scf(G^{n-1})\cong \underbrace{\scf(G)\otimes \cdots \otimes \scf(G)}_{\text{$n-1$ terms}}\otimes\ \scf(\{1\}).
\]

\subsection{Integer and set compositions} \label{Compositions}

An \emph{integer composition} $\mu=(\mu_1,\ldots,\mu_\ell)$ of $n$ is a sequence of positive integers $\mu_1,\ldots,\mu_\ell\in \ZZ_{\geq 1}$ such that $|\mu|=\mu_1+\cdots+\mu_\ell=n$.  We also write $\mu\vDash n$ and $\ell(\mu)=\ell$.  A convenient way to notate integer compositions is to give its ``ribbon" shape so as a sequence of rows of boxes appended end to end.  For example,
\[
(2,1,5,6,3,2)\longleftrightarrow 
\begin{tikzpicture}[scale=.2,baseline=.5cm]
\foreach \x/\y in {0/5,1/5,1/4,1/3,2/3,3/3,4/3,5/3,5/2,6/2,7/2,8/2,9/2,10/2,10/1,11/1,12/1,12/0,13/0}
	\draw (\x,\y) +(-.5,-.5) rectangle ++(.5,.5);
\end{tikzpicture}\ .
\]

It is well-known that the set of integer compositions of $n$ is in bijection with $\{0,1\}^{n-1}$.    However, there are two standard choices by reversing the roles of $0$ and $1$:
\[
\begin{array}{ccccc}   \{0,1\}^{n-1} & \longleftarrow &  \{\mu\vDash n\} & \longrightarrow  & \{0,1\}^{n-1}\\
 \binaryone(\mu) &  \mapsfrom & \mu & \mapsto & \binaryzero(\mu),
\end{array}
\]
where for $1\leq j\leq \ell(\mu)$,
\begin{align*}
 \binaryone(\mu)_j &= \left\{\begin{array}{ll}
 1 & \text{if $j\notin \{\mu_1,\mu_1+\mu_2,\ldots, \mu_1+\cdots + \mu_{\ell-1}\}$,}\\
  0 & \text{if $j\in \{\mu_1,\mu_1+\mu_2,\ldots, \mu_1+\cdots + \mu_{\ell-1}\}$,}\\
 \end{array}\right.\\ 
  \binaryzero(\mu)_j &= \left\{\begin{array}{ll}
 0 & \text{if $j\notin \{\mu_1,\mu_1+\mu_2,\ldots, \mu_1+\cdots + \mu_{\ell-1}\}$,}\\
  1 & \text{if $j\in \{\mu_1,\mu_1+\mu_2,\ldots, \mu_1+\cdots + \mu_{\ell-1}\}$.}\\
 \end{array}\right.
\end{align*}
In terms of ribbons, the binary sequence gives a sequence of right and down steps, so, for example,
\[
\begin{tikzpicture}[scale=.3,baseline=1cm]
\foreach \x/\h in {1/6,3/6,5/6,5/4,5/2,7/2,7/0,9/0,11/0,13/0,15/0}
	\draw (\x,\h) +(-.5,-.5) rectangle ++(.5,.5);
\foreach \x/\h/\z in {2/6/1,4/6/1,5/5/0,5/3/0,6/2/1,7/1/0,8/0/1,10/0/1,12/0/1,14/0/1}
	\node[gray] at (\x,\h) {$\scscs \z$};
\end{tikzpicture}
\overset{\binaryone}{\mapsfrom} \begin{tikzpicture}[scale=.3,baseline=.5cm]
\foreach \x/\h in {1/3,2/3,3/3,3/2,3/1,4/1,4/0,5/0,6/0,7/0,8/0}
	\draw (\x,\h) +(-.5,-.5) rectangle ++(.5,.5);
\end{tikzpicture} \overset{\binaryzero}{\mapsto}
\begin{tikzpicture}[scale=.3,baseline=1cm]
\foreach \x/\h in {1/6,3/6,5/6,5/4,5/2,7/2,7/0,9/0,11/0,13/0,15/0}
	\draw (\x,\h) +(-.5,-.5) rectangle ++(.5,.5);
\foreach \x/\h/\z in {2/6/0,4/6/0,5/5/1,5/3/1,6/2/0,7/1/1,8/0/0,10/0/0,12/0/0,14/0/0}
	\node[gray] at (\x,\h) {$\scscs \z$};
\end{tikzpicture}\ .
\]
In fact, if $\mu'$ denotes the conjugate composition of $\mu$, then
\[\binaryone(\mu)=\binaryzero(\mu').\]

Given two compositions $\mu=(\mu_1,\ldots, \mu_k)$ and $\nu=(\nu_1,\ldots,\nu_\ell)$, we have two ways of combining them:
\begin{description}
\item[Concate.]  $\mu\dotvdot\nu=(\mu_1,\ldots, \mu_k,\nu_1,\ldots,\nu_\ell)$, 
\item[Smash.] $\mu\vdotdotv\nu=(\mu_1,\ldots, \mu_{k-1},\mu_k+\nu_1,\nu_2,\ldots,\nu_\ell)$,
\end{description}
which correspond to adding a 0 or 1 between the corresponding binary sequences.

Using either of the above bijections in the case that $\scf(G)=\CC\spanning\{\mathds{1},\reg-\mathds{1}\}$, we obtain
\[
\dim(\scf(G^{n-1}))=\#\{\mu\vDash n\}.
\]

A \emph{set composition} $\underline{A}=(A_1,\ldots,A_\ell)$ of a set $B$ is a sequence of nonempty, pairwise disjoint subsets $A_1,\ldots, A_\ell\subseteq B$ such that 
\[
B=A_1\cup\cdots \cup A_\ell.
\]
We write $\underline{A}\vDash B$ and $\ell(\underline{A})=\ell$.  Fix a set $C$ and define a poset on the set compositions of $C$ by
\[
\underline{A}\preceq \underline{B}
\]
if for each $1\leq j\leq \ell(\underline{B})$ there exists an interval $1\leq h<i\leq \ell(\underline{A})$ such that $(A_h,A_{h+1},\ldots, A_i)\vDash B_j$.  For example,
\[
\begin{tikzpicture}
\node at (0,1) {$(\{3,5,6\},\{1,2,9,10,11\},\{4,7,8,12\})$};
\node at (0,0) {$(\{3,5,6\},\{9,10\},\{2,11\},\{1\},\{7,8\},\{4,12\})$};
\fill[gray!50] (-2.8,.8) -- (-3.2,.2) -- (-2.2,.2) -- (-1.8,.8) -- (-2.8,.8);
\fill[gray!50]  (-1.4,.8) -- (-1.8,.2) -- (1,.2) -- (.8,.8) -- (-1.4,.8);
\fill[gray!50] (1.2,.8) -- (1.4,.2)-- (3.2,.2) -- (2.8,.8) -- (1.2,.8);
\end{tikzpicture}
\]
In this case, we say that $\underline{A}$ \emph{refines} $\underline{B}$.  Another way to interpret the relation is that $\underline{A}$ is a sequence of set compositions on the blocks of $\underline{B}$.

\subsection{Combinatorial Hopf algebras}

The paper \cite{ABS} introduced the formal definition of a \emph{combinatorial Hopf algebra} as a graded connected Hopf algebra $\cH$ with an algebra homomorphism $\zeta:\cH\rightarrow \CC$, where $\CC$ is the underlying field of $\cH$.  While it is standard to refer to these homomorphisms as characters \cite{ABS}, due to the plethora of other characters in this paper we will refer to them as \emph{linear characters}.  The formal definition in fact makes every connected graded Hopf algebra a potential combinatorial one (by choosing $\zeta$ to be the counit of $\cH$).  The choice of $\zeta$ is combinatorial significant and each choice comes with a canonical Hopf algebra homomorphism to the Hopf algebra quasi-symmetric functions $\QSym$.  

The set of linear characters of $\cH$ forms a group under the convolution product 
\[
(\zeta\circ \xi)(h)=(m\circ(\zeta\otimes \xi)\circ \Delta)(h)\quad \text{for $h\in\cH$, $\zeta,\xi$ characters},
\]
but aside from this structure there is little guidance on what makes a ``good" choice of linear character. One of the benefits of the representation theoretic approach is that we easily obtain some natural choices along the way (we call ``supported by modules").  In Section \ref{CharacterGroupSection}, we even get a representation theoretic characterization of which functions $\zeta:\cH\rightarrow \CC$ are in fact linear characters (at least for our main example). 

As an element of a group the linear character $\zeta:\cH\rightarrow \CC$ also has an inverse, where the co-unit is the identity element of the group.  We say a character $\zeta$ is \emph{odd} if 
\[
\Res_{\cH_n}(\zeta^{-1})=(-1)^{n}\Res_{\cH_n}(\zeta),
\]
where $\cH_n$ is the $n$ degree component of $\cH$.

The main (well-known) examples of combinatorial Hopf algebras for this paper are as follows, where, for now, we omit any choice of corresponding linear character.

\vspace{.5cm}

\noindent\textbf{Noncommutative symmetric functions ($\NSym$).}  The Hopf algebra $\NSym$ is a graded non-commutative, cocommutative Hopf algebra with bases
\[\NSym=\CC\spanning\{H_\mu\mid \mu\vDash n,n\in \ZZ_{\geq 0}\}=\CC\spanning\{R_\mu\mid \mu\vDash n, n\in \ZZ_{\geq 0}\},\]
where the latter is known as the ribbon basis.  The Hopf structure is given by
\[
H_\mu\cdot H_\nu=H_{\mu\sdotvdot\nu}, \quad \text{where $\mu\dotvdot\nu$ is concatenation as in Section \ref{Compositions},}
\]
and 
\[
\Delta(H_{(n)})=\sum_{j=0}^n H_j\otimes H_{n-j}.
\]
The second basis is given by
\[
H_\mu=\sum_{\nu \text{ coarsens } \mu} R_\nu,
\]
with the product
\[R_\mu \cdot R_\nu=R_{\mu\sdotvdot \nu}+R_{\mu\svdotdotv \nu}.\]
One goal of this paper is to give an interpretation of this Hopf algebra via the representation theory of a tower of groups.  The dual Hopf algebra $\QSym=\NSym^\dual$ is the Hopf algebra of \emph{quasi-symmetric functions}.  

\vspace{.5cm}

\noindent\textbf{Symmetric functions in noncommuting variables ($\NCSym$).} The Hopf algebra $\NCSym$ is a graded non-commutative, cocommutative Hopf algebra with basis
\[
\NCSym=\CC\spanning\{m_\mu\mid \mu \text{ a set partition of $\{1,2,\ldots, n\}$, $n\in \ZZ_{\geq 0}$}\},
\]
and can be thought of as a version of symmetric functions $\Sym$ where the variables no longer commute.  The paper \cite{Aim} shows that this Hopf algebra can be realized representation theoretically by considering supercharacter theories 
\[
\scf(\UT_\bullet)=\bigoplus_{n\geq 0} \scf(\UT_n(\FF_q)),
\]
where $\UT_n(\FF_q)$ is the finite group of $n\times n$ unipotent  upper triangular matrices over the field $\FF_q$ with $q$ elements.  In this case, the basis $\{m_\lambda\}_{|\lambda|=n}$ corresponds to the basis of superclass identifier functions of $\scf(\UT_n(\FF_q))$.

\vspace{.5cm}

\noindent\textbf{Malvenuto--Reutenauer Hopf algebra ($\MR$).} The Malvenuto--Reutenauer algebra $\MR$ is a noncommutative, noncocommutive self dual Hopf algebra with basis
\[
\MR=\CC\spanning\{F_w\mid w\in S_n, n\in \ZZ_{\geq 0}\}. 
\]
The paper \cite{ATmr} shows that this Hopf algebra can be realized representation theoretically by considering supercharacter theories 
\[
\scf(\fkut_\bullet)=\bigoplus_{n\geq 0} \scf(\fkut_n(\FF_q)),
\]
where $\fkut_n(\FF_q)$ is the additive group of $n\times n$ nilpotent upper triangular matrices.  Here, the supercharacters of $\fkut_n(\FF_q)$ correspond to the fundamental basis $\{F_w\}_{w\in S_n}$.

\section{Functorial constructions for Hopf structures}\label{HopfAlgebraConstruction}

Fix a finite group $G$ with a supercharacter theory $(\Ch,\Cl)$ and corresponding space of superclass functions $\scf(G)$.  If we let  $\scf(G^{n-1})$ comes from the supercharacter theory defined in (\ref{DirectProductSupercharacterTheory}), we obtain a vector space
\[
\scf(G^\all)=\bigoplus_{n\geq 0} \scf(G^{n-1}),
\]  
where by convention
\[
 \scf(G^{1-1})=\scf(\{1\})\quad \text{and}\quad \scf(G^{0-1})=\CC\spanning\{\chi^{\emptyset}\},
\]
and $\chi^\emptyset$ will be the unit of our Hopf algebra below.  The following section imposes a Hopf structure on this vector space using representation theoretic functors.  As a running example, we will use the case $\Ch=\{\mathds{1},\reg-\mathds{1}\}$.

\begin{remark}
At first glance, the grading we have chosen seems a bit odd, since it might be more natural to let the $n$th graded component be $\scf(G^n)$.  However, this choice gives far less flexibility representation theoretically since $G^m\times G^n=G^{m+n}$, where-as $G^{m-1}\times G^{n-1}\subseteq G^{m+n-1}$ (with the containment generally strict).  We were initially motivated by the abelianization $\UT_n(\FF_q)$, given by
\[
\UT_n(\FF_q)/[\UT_n(\FF_q),\UT_n(\FF_q)]\cong \underbrace{\FF_q^+\times \cdots \times \FF_q^+}_{\text{$n-1$ terms}}.
\]
In this case, we have that the dotted arrow
\[
\begin{tikzpicture}[baseline=1.5cm]
\node (02) at (0,3) {$\UT_{m+n}$};
\node (11)  at (-3,1.5) {$\UT_m\times \UT_n$};
\node (01) at (3,1.5) {$\UT_{m+n}/[\UT_{m+n},\UT_{m+n}]$};
\node (00) at (0,0) {$\UT_{m}/[\UT_{m},\UT_{m}]\times \UT_{n}/[\UT_{n},\UT_{n}]$};
\draw[->] (11) -- (02);
\draw[->] (02) -- (01);
\draw[->] (11) -- (00);
\draw[->,dotted] (00) -- (01);
\end{tikzpicture},
\]
making this diagram commute is exactly our inclusion from our chosen grading.
\end{remark}

\subsection{Basic functors}

All our functors will be based on the interaction between a set of normal subgroups of $G^{n-1}$ in bijection with $\{0,1\}^{n-1}$.  Specifically, if $\underline{a}=(a_1,\ldots, a_{n-1})\in \{0,1\}^{n-1}$, let
\[
G^{n-1}_{\underline{a}}=\{(g_1,\ldots,g_{n-1},1)\in G^{n-1}\mid g_j\neq 1 \text{ implies } a_j=1\},
\]
so that the lattice of subgroups is exactly isomorphic to the Boolean lattice on $\{0,1\}^{n-1}$.    For example,
\[
G^{10-1}_{(1,0,0,1,1,0,1,0,1)}=G\times \{1\}\times\{1\} \times G\times G\times \{1\} \times G\times \{1\}\times G\times \{1\}.
\]

Not only is $G^{n-1}_{\underline{a}}\subseteq G^{n-1}$ a subgroup, but since $G^{n-1}_{\underline{a}}$ has a normal complement $G^{n-1}_{\overline{\underline{a}}}$ ($\overline{\underline{a}}$ is the binary sequence obtained by switching zeroes and ones), it is also a quotient under a corresponding canonical projection $\pi_{\underline{a}}$.  We therefore have standard functors induction and inflation,
\[
\Ind:G_{\underline{a}}^{n-1}\text{-mod} \longrightarrow G^{n-1}\text{-mod}\quad \text{and}\quad \Inf:G_{\underline{a}}^{n-1}\text{-mod} \longrightarrow G^{n-1}\text{-mod}
\]
with corresponding Frobenius duals restriction $\Res$ and deflation $\Def$, respectively.  The corresponding transformations on superclass functions are given by
\[
\begin{array}{rccc} \Ind_{G_{\underline{a}}^{n-1}}^{G^{n-1}}:  & \scf(G_{\underline{a}}^{n-1}) & \longrightarrow & \scf(G^{n-1})\\
& \psi & \mapsto & \underline{g}\mapsto \left\{\begin{array}{ll}\frac{|G^{n-1}|}{|G^{n-1}_{\underline{a}}|} \psi(\underline{g}) & \text{if $\underline{g}\in G_{\underline{a}}^{n-1}$}\\ 0 & \text{otherwise.}\end{array}\right.\end{array}
\]
\[
\begin{array}{rccc} \Inf_{G_{\underline{a}}^{n-1}}^{G^{n-1}}:  & \scf(G_{\underline{a}}^{n-1}) & \longrightarrow & \scf(G^{n-1})\\
& \psi & \mapsto & \psi\circ \pi_{\underline{a}}.
\end{array}
\]
By inspection, we obtain the following result, where we use $\odot$ to indicate pointwise multiplication (corresponding to the diagonal action on modules).
\begin{lemma} For $\underline{a}\in \{0,1\}^{n-1}$,
\[
 \Ind_{G_{\underline{a}}^{n-1}}^{G^{n-1}}(\cdot)= \Inf_{G_{\overline{\underline{a}}}^{n-1}}^{G^{n-1}}(\reg) \odot \Inf_{G_{\underline{a}}^{n-1}}^{G^{n-1}}(\cdot),
\]
where $\reg$ is the regular character of $G_{\overline{\underline{a}}}^{n-1}$.
\end{lemma} 
By Frobenius reciprocity, we obtain the complementary result.
\begin{corollary}
 For $\underline{a}\in \{0,1\}^{n-1}$,
\[
 \Res_{G_{\underline{a}}^{n-1}}^{G^{n-1}}(\cdot)= \Def_{G_{\underline{a}}^{n-1}}^{G^{n-1}}\Big(\Inf_{G_{\overline{\underline{a}}}^{n-1}}^{G^{n-1}}(\reg) \odot \cdot \Big)
\]
\end{corollary}

Thus, if we wish to work in a setting where we include all the standard functors, we can consider functions of the form
\[
\psi\mapsto \alpha\odot \Inf_{G_{\overline{\underline{a}}}^{n-1}}^{G^{n-1}}(\psi) \quad \text{and}\quad \psi\mapsto  \Def_{G_{\overline{\underline{a}}}^{n-1}}^{G^{n-1}}(\beta\odot \psi)
\]
for $\alpha,\beta\in \scf(G^{n-1})$.

We say a function $\psi\in \scf(G^{n-1})$ \emph{factors} if there are functions $\psi_1,\psi_2,\ldots,\psi_{n-1}\in \scf(G)$ such that 
\[
\psi(g_1,\ldots, g_{n-1},1)=\psi_1(g_1)\psi_2(g_2)\cdots \psi_{n-1}(g_{n-1}) \quad \text{for $g_1,\ldots,g_{n-1}\in G$.}
\]
Equivalently, $\psi$ is a simple tensor in the isomorphic space $\scf(G)\otimes \cdots \otimes \scf(G)\otimes \scf(\{1\})$, so we will write $\psi=\psi_1\otimes \psi_2\otimes\cdots\times \otimes \psi_{n}$ (where $\psi_{n}=\chi^{()}$).  Note that the elements in our two canonical bases (supercharacters and superclass identifier functions) both factor, by the construction of the supercharacter theory.  

As a group,
\[
G^{n-1}_{\underline{a}}\cong G^{|\underline{a}|},\quad \text{where}\quad |\underline{a}|=a_1+a_2+\cdots+a_{n-1}.
\]
We can therefore use the same supercharacter theories for these groups.  However, as a subgroup of $G^{n-1}$ it is often convenient to factor functions $\gamma =\gamma_1\otimes \cdots \otimes \gamma_{n}\in \scf(G_{\underline{a}}^{n-1})$ into $n$ terms, where $\gamma_j=\chi^{()}$ if $a_j=0$ or $j=n$.

\begin{proposition} Let $\underline{a}\in \{0,1\}^{n-1}$ with $|\underline{a}|=m$.  Then
\begin{enumerate}
\item[(a)] If $\gamma=\gamma_1\otimes \gamma_2\otimes \cdots \otimes \gamma_{n}\in \scf(G^{n-1}_{\underline{a}})$ factors, then so does $\Inf_{G_{\overline{\underline{a}}}^{n-1}}^{G^{n-1}}(\gamma)$.
\item[(b)] Given a supercharacter $\chi=\chi_1\otimes \chi_2\otimes \cdots \otimes \chi_{n}\in \scf(G^{n-1}_{\underline{a}})$ and $1\leq j\leq n-1$,
\[\Big(\Inf_{G_{\overline{\underline{a}}}^{n-1}}^{G^{n-1}}(\chi)\Big)_j=\left\{\begin{array}{ll}
\chi_j & \text{if $a_j=1$},\\
\mathds{1} & \text{if $a_j=0$.}
\end{array}\right.
\]
\item[(c)] Given a superclass identifier function $\delta=\delta_1\otimes \delta_2\otimes \cdots \otimes \delta_{n}\in \scf(G^{n-1}_{\underline{a}})$ and $1\leq j\leq n-1$,
\[\Big(\Inf_{G_{\overline{\underline{a}}}^{n-1}}^{G^{n-1}}(\delta)\Big)_j=\left\{\begin{array}{ll}
\delta_j & \text{if $a_j=1$},\\
\sum_{S\in \Cl} \delta_S & \text{if $a_j=0$.}
\end{array}\right.
\]
\end{enumerate}
\end{proposition}

We have a similar result for deflation,
\begin{proposition} Let $\underline{a}\in \{0,1\}^{n-1}$ with $|\underline{a}|=m$.  Then
\begin{enumerate}
\item[(a)] If $\gamma=\gamma_1\otimes \gamma_2\otimes \cdots \otimes \gamma_{n}\in \scf(G^{n-1})$ factors, then so does $\Def_{G_{\overline{\underline{a}}}^{n-1}}^{G^{n-1}}(\gamma)$.
\item[(b)] Given a supercharacter $\chi=\chi_1\otimes \chi_2\otimes \cdots \otimes \chi_{n}\in \scf(G^{n-1})$ and $1\leq j\leq n-1$,
\[\Big(\Def_{G_{\overline{\underline{a}}}^{n-1}}^{G^{n-1}}(\chi)\Big)_j=\left\{\begin{array}{ll}
\chi_j & \text{if $a_j=1$},\\
\langle\chi_j,\mathds{1}\rangle \chi^{()} & \text{if $a_j=0$.}
\end{array}\right.
\]
\item[(c)] Given a superclass identifier function $\delta=\delta_1\otimes \delta_2\otimes \cdots \otimes \delta_{n}\in \scf(G^{n-1}_{\underline{a}})$ and $1\leq j\leq n-1$,
\[\Big(\Def_{G_{\overline{\underline{a}}}^{n-1}}^{G^{n-1}}(\delta)\Big)_j=\left\{\begin{array}{ll}
\delta_j & \text{if $a_j=1$},\\
\frac{|S_j|}{|G|}\chi^{()} & \text{if $a_j=0$,}
\end{array}\right.
\]
where $S_j$ is the superclass identified by $\delta_j$.
\end{enumerate}
\end{proposition}

\subsection{Product functor}

Given a set composition $\underline{A}=(A_1,A_2,\ldots,A_\ell)$ of $\{1,2,\ldots, n\}$, let $\lc(\underline{A})\in \{0,1\}^{n-1}$ be the binary sequence given by
\[
\lc(\underline{A})_j=\left\{\begin{array}{ll}
1 & \text{if $j \notin\{\max(A_1),\max(A_2),\ldots,\max(A_\ell)\}$,}\\
0 & \text{otherwise.}
\end{array}\right.
\]
For example, if $\underline{A}=(\{1,3,4,5,7\},\{6\},\{8,9\},\{2,10\})$, then
\[
\begin{tikzpicture}[scale=.5,baseline=-.5cm]
\foreach \x in {1,...,9}
	\node at (\x+.5,0) {$\cdot$};
\foreach \x in {1,...,10}
	\node at (\x,0) {$\shortmid$};
\foreach \x in {1,3,4,5,7}
	\node (\x) at (\x,-.5) [inner sep=1pt] {$\scs\x$};
\foreach \x in {6}
	\node  (\x) at (\x,-1) {$\scs\x$};
\foreach \x in {8,9}
	\node (\x) at (\x,-1.5)  [inner sep=1pt] {$\scs\x$};
\foreach \x in {2,10}
	\node  (\x) at (\x,-2)  [inner sep=1pt] {$\scs\x$};
\foreach \x in {1,3,4,5,8,2}
	\node at (\x+.5,.5)  [inner sep=1pt] {$\scs 1$};
\foreach \x in {7,6,9}
	\node at (\x+.5,.5)  {$\scs 0$};
\foreach \x in {1,2,3,4,5,8}
	\draw[gray,very thick] (\x) -- ++(.8,0);
\node at (0,.5) {$\scs\lc(\underline{A})$};
\node at (0,-.5) {$\scs A_1$};
\node at (0,-1) {$\scs A_2$};
\node at (0,-1.5) {$\scs A_3$};
\node at (0,-2) {$\scs A_4$};
\end{tikzpicture}\ .
\]
Thus we view the nonmaximal block elements as selecting the left coordinates of the entries that will be 1.   In this case,
\[
G^{10-1}_{\lc(\underline{A})} = G\times G\times G\times G\times G\times \{1\}\times  \{1\} \times G\times \{1\} \cong G^{5-1}\times G^{1-1}\times G^{2-1}\times G^{2-1},
\]
where the isomorphism comes from rearranging the factors block by block.

For $\iota\in \scf(G)$, $\underline{A},\underline{B}\vDash\{1,2,\ldots,n\}$ with $\underline{A}$ a refinement of $\underline{B}$, define
\[
\begin{array}{rccc} \Inf_{\underline{A}}^{\underline{B}}[\iota]: & \scf(G^{n-1}_{\lc(\underline{A})}) & \longrightarrow & \scf(G^{n-1}_{\lc(\underline{B})})\\
&\psi & \mapsto & \iota_{\underline{A}}^{\underline{B}}\odot \Inf_{G^{n-1}_{\lc(\underline{A})}}^{G^{n-1}_{\lc(\underline{B})}}(\psi),
\end{array}
\]
where $\iota_{\underline{A}}^{\underline{B}}\in \scf(G^{n-1}_{\lc(\underline{B})})$ is a factorizable function given by
\[
(\iota_{\underline{A}}^{\underline{B}})_j=\left\{\begin{array}{ll}
\mathds{1} & \text{if $\lc(\underline{A})_j=\lc(\underline{B})_j=1,$}\\
\iota & \text{if $\lc(\underline{A})_j<\lc(\underline{B})_j$,}\\
\chi^{()} & \text{if $\lc(\underline{A})_j=\lc(\underline{B})_j=0$.}
\end{array}\right.
\]
The function $\iota_{\underline{A}}^{\underline{B}}$ effectively inserts an $\iota$ in each coordinate $j$ where $\lc(\underline{A})_j<\lc(\underline{B})_j$.
Note that $\Inf_{\underline{A}}^{\underline{B}}[\mathds{1}]$ is ordinary inflation, and  $\Inf_{\underline{A}}^{\underline{B}}[\reg]$ is induction.

These functions are transitive in the ways one would hope.
\begin{lemma}\label{LemmaAssociativity}
Let $\underline{A},\underline{B}\vDash \{1,2,\ldots,n\}$ such that  $\underline{A}$ refines $\underline{B}$.  
\begin{enumerate}
\item[(a)]  If $\psi\in \scf(G^{n-1}_{\lc(\underline{A})})$ factors, then $\Inf_{\underline{A}}^{\underline{B}}[\iota] (\psi)$ factors.
\item[(b)]  If $\underline{C}\vDash\{1,2,\ldots, n\}$ such that $\underline{B}$ refines $\underline{C}$, then
 \[\Inf_{\underline{B}}^{\underline{C}}[\iota] \circ\Inf_{\underline{A}}^{\underline{B}}[\iota] =\Inf_{\underline{A}}^{\underline{C}}[\iota].
\]
\end{enumerate}
\end{lemma}
\begin{proof}
(a) Suppose $\psi=\psi_1\otimes\cdots\otimes \psi_n$.  Then
\[
\Inf_{\underline{A}}^{\underline{B}}[\iota] (\psi)=\Big((\iota_{\underline{A}}^{\underline{B}})_1\odot \Inf_{G_{\lc(\underline{A})_1}^{n-1}}^{G_{\lc(\underline{B})_1}^{n-1}}(\psi_1)\Big)\otimes \Big((\iota_{\underline{A}}^{\underline{B}})_2\odot \Inf_{G_{\lc(\underline{A})_2}^{n-1}}^{G_{\lc(\underline{B})_2}^{n-1}}(\psi_2)\Big)\otimes \cdots \otimes \Big((\iota_{\underline{A}}^{\underline{B}})_n\odot \Inf_{G_{\lc(\underline{A})_n}^{n-1}}^{G_{\lc(\underline{B})_n}^{n-1}}(\psi_n)\Big).
\]
(b) If $\psi=\psi_1\otimes\cdots\otimes \psi_n$,
\begin{align*}
\Big(\Inf_{\underline{B}}^{\underline{C}}[\iota] \circ\Inf_{\underline{A}}^{\underline{B}}[\iota](\psi)\Big)_j 
&=(\iota_{\underline{B}}^{\underline{C}})_j\odot \Inf_{G_{\lc(\underline{B})_j}^{n-1}}^{G_{\lc(\underline{C})_j}^{n-1}} \Big((\iota_{\underline{A}}^{\underline{B}})_j\odot \Inf_{G_{\lc(\underline{A})_j}^{n-1}}^{G_{\lc(\underline{B})_j}^{n-1}} (\psi_j)\Big)\\
&=\left\{\begin{array}{ll}
\mathds{1}\odot\mathds{1}\odot \psi_j & \text{if $\lc(\underline{A})_j=\lc(\underline{B})_j=\lc(\underline{C})_j=1$}\\
\mathds{1}\odot\iota\odot \mathds{1} & \text{if $\lc(\underline{A})_j<\lc(\underline{B})_j=\lc(\underline{C})_j=1$}\\
\iota\odot \mathds{1} & \text{if $\lc(\underline{A})_j=\lc(\underline{B})_j<\lc(\underline{C})_j=1$}\\
\chi^{()} \odot\chi^{()} \odot \chi^{()} & \text{if $\lc(\underline{A})_j=\lc(\underline{B})_j=\lc(\underline{C})_j=0$}\\
\end{array}\right.\\
&=\Inf_{\underline{A}}^{\underline{C}}[\iota](\psi)_j.
\end{align*}
\end{proof}

We define a product on $\scf(G^\all)$ by extending the graded product
\begin{equation}\label{HopfProduct}
\begin{array}{ccc} \scf(G^{m-1})\otimes \scf(G^{n-1}) & \longrightarrow & \scf(G^{m+n-1})\\
\gamma \otimes \psi & \mapsto & \Inf_{(\{1,2,\ldots, m\},\{m+1,\ldots,n\})}^{(\{1,2,\ldots,m+n\})} [\iota](\gamma\otimes \psi)
\end{array}
\end{equation}
linearly.  By Lemma \ref{LemmaAssociativity} (b), this product is associative.  

Note that if $\psi=\psi_1\otimes \cdots\otimes \psi_m$ and $\gamma=\gamma_1\otimes \cdots\otimes \gamma_n$ factor, then
\[
\Inf_{(\{1,2,\ldots, m\},\{m+1,\ldots,n\})}^{(\{1,2,\ldots,m+n\})} [\iota](\psi\otimes \gamma)=\psi_1\otimes \cdots\otimes \psi_{m-1}\otimes \iota \otimes \gamma_1\otimes \cdots\otimes \gamma_n.
\]
If we want to work in a specific basis, therefore, we merely need to know how to decompose $\iota$.

\begin{runningexample}
In our example with $\scf(G)=\CC\spanning\{\mathds{1},\reg\}$, 
\begin{description}
\item[Case $\iota=\mathds{1}$.]  For $\underline{a}\in \{0,1\}^{m-1}$ and $\underline{b}\in \{0,1\}^{n-1}$,
\begin{align*}
\chi^{\underline{a}}\cdot \chi^{\underline{b}}&=\chi^{({\underline{a}},1,{\underline{b}})}\\
\delta_{\underline{a}}\cdot  \delta_{\underline{b}} &= \delta_{({\underline{a}},1,{\underline{b}})}+\delta_{({\underline{a}},0,{\underline{b}})}.
\end{align*}
\item[Case $\iota=\reg$.]  For $\mu\vDash m$ and $\nu\vDash n$,
\begin{align*}
\chi^{\underline{a}}\cdot \chi^{\underline{b}}&=\chi^{({\underline{a}},0,{\underline{b}})}+\chi^{({\underline{a}},1,{\underline{b}})}\\
\delta_{\underline{a}}\cdot  \delta_{\underline{b}} &=|G| \delta_{({\underline{a}},0,{\underline{b}})}.
\end{align*}
\end{description}
There is in fact a Pontryagin style duality between the two cases which encodes why these results seem so similar (while being reversed).
\end{runningexample}

\subsection{Coproduct functor}
Given a set composition $\underline{A}=(A_1,A_2,\ldots, A_\ell)\vDash\{1,2,\ldots, n\}$,  define binary sequences $\llc(\underline{A}),\bc(\underline{A})\in \{0,1\}^{n-1}$ by
\begin{align*}
\llc(\underline{A})_j & =\left\{\begin{array}{ll} 1 & \text{if $j\in A_k$, then $j+1\in A_l$ with $k\leq l$,}\\
0 & \text{otherwise,}\end{array}\right.\\
\bc(\underline{A})_j & =\left\{\begin{array}{ll} 1 & \text{if $j\in A_k$, then $j+1\in A_k$,}\\
0 & \text{otherwise.}\end{array}\right.
\end{align*}
For example, if $\underline{A}=(\{1,3,4,5,7\},\{6\},\{8,9\},\{2,10\})$, then
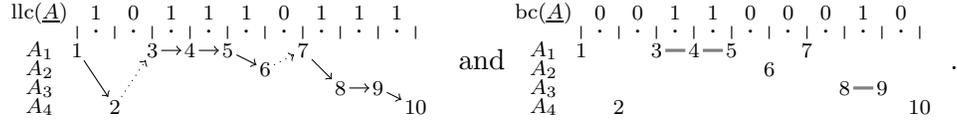
\[
\begin{tikzpicture}[scale=.5,baseline=-.5cm]
\foreach \x in {1,...,9}
	\node at (\x+.5,0) {$\cdot$};
\foreach \x in {1,...,10}
	\node at (\x,0) {$\shortmid$};
\foreach \x in {1,3,4,5,7}
	\node (\x) at (\x,-.5) [inner sep=1pt] {$\scs\x$};
\foreach \x in {6}
	\node (\x) at (\x,-1) [inner sep=1pt] {$\scs\x$};
\foreach \x in {8,9}
	\node (\x) at (\x,-1.5) [inner sep=1pt] {$\scs\x$};
\foreach \x in {2,10}
	\node (\x) at (\x,-2) [inner sep=1pt] {$\scs\x$};
\foreach \x in {1,3,4,5,8,7,9}
	\node at (\x+.5,.5) {$\scs 1$};
\foreach \x in {2,6}
	\node at (\x+.5,.5) {$\scs 0$};
\node at (0,.5) {$\scs\llc(\underline{A})$};
\node at (0,-.5) {$\scs A_1$};
\node at (0,-1) {$\scs A_2$};
\node at (0,-1.5) {$\scs A_3$};
\node at (0,-2) {$\scs A_4$};
\foreach \x/\y in {1/2,3/4,4/5,5/6,7/8,8/9,9/10}
\draw[->] (\x) -- (\y); 
\foreach \x/\y in {2/3,6/7}
\draw[->, dotted] (\x) -- (\y); 
\end{tikzpicture}
\quad\text{and}
\begin{tikzpicture}[scale=.5,baseline=-.5cm]
\foreach \x in {1,...,9}
	\node at (\x+.5,0) {$\cdot$};
\foreach \x in {1,...,10}
	\node at (\x,0) {$\shortmid$};
\foreach \x in {1,3,4,5,7}
	\node (\x) at (\x,-.5) [inner sep=1pt]  {$\scs\x$};
\foreach \x in {6}
	\node at (\x,-1) {$\scs\x$};
\foreach \x in {8,9}
	\node (\x) at (\x,-1.5) [inner sep=1pt] {$\scs\x$};
\foreach \x in {2,10}
	\node at (\x,-2) {$\scs\x$};
\foreach \x in {3,4,8}
	\node at (\x+.5,.5) {$\scs 1$};
\foreach \x in {1,2,5, 6,7,9}
	\node at (\x+.5,.5) {$\scs 0$};
\node at (0,.5) {$\scs\bc(\underline{A})$};
\node at (0,-.5) {$\scs A_1$};
\node at (0,-1) {$\scs A_2$};
\node at (0,-1.5) {$\scs A_3$};
\node at (0,-2) {$\scs A_4$};
\foreach \x/\y in {3/4,4/5,8/9}
\draw[gray,very thick] (\x) -- (\y); 
\end{tikzpicture}\ .
\]
We view $\llc(\underline{A})$ as selecting those left coordinates $j$ for which $j+1$ is in the same block or a subsequent block, and $\bc(\underline{A})$ as only allowing those coordinates $j$ where $j+1$ is in the same block.

For $\tau, \alpha,\beta \in \scf(G)$, $\underline{A},\underline{B}\vDash\{1,2,\ldots,n\}$ with $\underline{A}$ a refinement of $\underline{B}$, define
\[
\begin{array}{rccc} \Down_{\underline{A}}^{\underline{B}}[\tau,\alpha,\beta]: & \scf(G^{n-1}_{\lc(\underline{B})}) & \longrightarrow & \scf(G^{n-1}_{\lc(\underline{A})})\\
&\psi & \mapsto & \tau_{\underline{A}}\odot \Inf_{G^{n-1}_{\bc(\underline{A})}}^{G^{n-1}_{\lc(\underline{A})}}\circ \Def_{G^{n-1}_{\bc(\underline{A})}}^{G^{n-1}_{\lc(\underline{B})}}((\alpha\times\beta)_{\underline{A}}^{\underline{B}}\odot\psi),
\end{array}
\]
where $\tau_{\underline{A}}\in \scf(G^{n-1}_{\lc(\underline{A})})$ and $(\alpha\times\beta)_{\underline{A}}^{\underline{B}}\in \scf(G^{n-1}_{\lc(\underline{B})})$ are the factoring functions given by
\[
(\tau_{\underline{A}})_j=\left\{\begin{array}{ll}
\mathds{1} & \text{if $\bc(\underline{A})_j=\lc(\underline{A})_j=1,$}\\
\tau & \text{if $\bc(\underline{A})_j<\lc(\underline{A})_j$},\\
\chi^{()} &  \text{if $\bc(\underline{A})_j=\lc(\underline{A})_j=0,$}\\
\end{array}\right.\quad\text{and}\quad
((\alpha\times \beta)_{\underline{A}}^{\underline{B}})_j=\left\{\begin{array}{ll}
\mathds{1} & \text{if $\bc(\underline{A})_j=1$,}\\
\alpha & \text{if $\bc(\underline{A})_j<\llc(\underline{A})_j$},\\
\beta & \text{if $0=\llc(\underline{A})_j<\lc(\underline{B})_j$},\\
\chi^{()} & \text{if $0=\llc(\underline{A})_j=\lc(\underline{B})_j$}.
\end{array}\right.
\]
Here we have $\Down[\mathds{1},\mathds{1},\mathds{1}]=\Inf\circ\Def$, $\Down[\mathds{1},\reg,\reg]=\Inf\circ\Res$, and $\Down[\reg,\mathds{1},\mathds{1}]=\Ind\circ\Def$. 

These functions are transitive in the ways one would hope as long as we add some conditions on $\tau$, $\alpha$ and $\beta$.
\begin{lemma}\label{LemmaCoAssociativity}
Let $\underline{A},\underline{B}\vDash \{1,2,\ldots,n\}$ such that  $\underline{A}$ refines $\underline{B}$.  
\begin{enumerate}
\item[(a)] If $\psi\in \scf(G^{n-1}_{\lc(\underline{B})})$ factors, then $\Down_{\underline{A}}^{\underline{B}}[\iota] (\psi)$ factors.
 \item[(b)] Suppose $\langle\tau,\alpha\rangle=\langle\tau,\beta\rangle=1$.  If $\underline{C}\vDash\{1,2,\ldots, n\}$ such that $\underline{B}$ refines $\underline{C}$, then
\[\Down_{\underline{A}}^{\underline{B}}[\tau,\alpha,\beta] \circ\Down_{\underline{B}}^{\underline{C}}[\tau,\alpha,\beta] =\Down_{\underline{A}}^{\underline{C}}[\tau,\alpha,\beta].
\]
\end{enumerate}
\end{lemma}
\begin{proof}
(a) Suppose $\psi=\psi_1\otimes\cdots\otimes \psi_n$.  Then
\[
\Down_{\underline{A}}^{\underline{B}}[\tau,\alpha,\beta] (\psi)
=\bigotimes_{j=1}^n\Big((\tau_{\underline{A}})_j\odot \Inf_{G^{n-1}_{\bc(\underline{A})_j}}^{G^{n-1}_{\lc(\underline{A})_j}}\circ \Def_{G^{n-1}_{\bc(\underline{A})_j}}^{G^{n-1}_{\lc(\underline{B})_j}}(((\alpha\times\beta)_{\underline{A}}^{\underline{B}})_j\odot\psi_j)\Big).
\]
(b) We begin with some observations about refinements.  First, if $\underline{A}$ refines $\underline{B}$, then $\bc(\underline{A})_j\leq \bc(\underline{B})_j$ and $\lc(\underline{A})_j\leq \lc(\underline{B})_j$, and $\llc(\underline{A})_j\leq \llc(\underline{B})_j$ for all $1\leq j\leq n$.  Second, that even stronger condition of $\bc(\underline{B})_j=0$ implies $\llc(\underline{A})_j=\llc(\underline{B})_j$ for all $1\leq j\leq n$ holds, since $\underline{A}$ refines $\underline{B}$ block by block.

If $\psi=\psi_1\otimes\cdots\otimes \psi_n$,  
\begin{align*}
&\Big(\Down_{\underline{A}}^{\underline{B}}[\tau,\alpha,\beta]\circ\Down_{\underline{B}}^{\underline{C}}[\tau,\alpha,\beta](\psi)\Big)_j\\ 
&=(\tau_{\underline{A}})_j\odot \Inf_{G^{n-1}_{\bc(\underline{A})_j}}^{G^{n-1}_{\lc(\underline{A})_j}}\hspace{-.2cm}\circ \Def_{G^{n-1}_{\bc(\underline{A})_j}}^{G^{n-1}_{\lc(\underline{B})_j}}\hspace{-.1cm}\Big(((\alpha\times\beta)_{\underline{A}}^{\underline{B}})_j\odot(\tau_{\underline{B}})_j\odot \Inf_{G^{n-1}_{\bc(\underline{B})_j}}^{G^{n-1}_{\lc(\underline{B})_j}}\hspace{-.2cm}\circ \Def_{G^{n-1}_{\bc(\underline{B})_j}}^{G^{n-1}_{\lc(\underline{C})_j}}\hspace{-.1cm}(((\alpha\times\beta)_{\underline{B}}^{\underline{C}})_j\odot\psi_j)\Big)\\
&=\left\{\begin{array}{ll}
\psi_j & \text{if $\bc(\underline{A})_j=1$,}\\
\langle\psi_j,\alpha\rangle \tau & \text{if $\bc(\underline{A})_j<\bc(\underline{B})_j$, $\llc(\underline{A})_j=\lc(\underline{A})_j=1$,}\\
\langle\psi_j,\alpha\rangle \chi^{()} & \text{if $\bc(\underline{A})_j<\bc(\underline{B})_j$, $\llc(\underline{A})_j>\lc(\underline{A})_j$,}\\
\langle\psi_j,\beta\rangle \tau & \text{if $\bc(\underline{A})_j<\bc(\underline{B})_j$, $\llc(\underline{A})_j<\lc(\underline{A})_j$,}\\
\langle\psi_j,\beta\rangle\chi^{()} & \text{if $\bc(\underline{A})_j<\bc(\underline{B})_j$, $\llc(\underline{A})_j=\lc(\underline{A})_j=0$,}\\
\langle\psi_j,\alpha\rangle \langle\tau,\alpha\rangle\tau & \text{if $\bc(\underline{B})_j<\llc(\underline{A})_j=\lc(\underline{A})_j$,}\\
\langle\psi_j,\alpha\rangle \langle\tau,\alpha\rangle\chi^{()} & \text{if $\bc(\underline{B})_j=\lc(\underline{A})_j<\llc(\underline{A})_j=\lc(\underline{B})_j$,}\\
\langle\psi_j,\alpha\rangle \chi^{()} & \text{if $\bc(\underline{B})_j=\lc(\underline{B})_j<\llc(\underline{A})_j$,}\\
\langle\psi_j,\beta\rangle \langle\tau,\beta\rangle\tau & \text{if $\bc(\underline{B})_j=\llc(\underline{B})_j<\lc(\underline{A})_j$,}\\
\langle\psi_j,\beta\rangle \langle\tau,\beta\rangle\chi^{()}& \text{if $\bc(\underline{B})_j=\llc(\underline{B})_j=\lc(\underline{A})_j<\lc(\underline{B})_j$,}\\
\langle\psi_j,\beta\rangle\chi^{()}& \text{if $\bc(\underline{B})_j=\llc(\underline{B})_j=\lc(\underline{B})_j<\lc(\underline{C})_j$,}\\
\chi^{()}& \text{if $\lc(\underline{C})_j=0$,}\\
\end{array}\right.\\
&=(\tau_{\underline{A}})_j\odot \Inf_{G^{n-1}_{\bc(\underline{A})_j}}^{G^{n-1}_{\lc(\underline{A})_j}}\circ \Def_{G^{n-1}_{\bc(\underline{A})_j}}^{G^{n-1}_{\lc(\underline{C})_j}}(((\alpha\times\beta)_{\underline{A}}^{\underline{C}})_j\odot\psi_j).
\end{align*}
\end{proof}

Define a coproduct on $\scf(G^\all)$ by extending
\begin{equation}\label{HopfCoProduct}
\begin{array}{ccc}  \scf(G^{n-1}) & \longrightarrow & \dd\bigoplus_{m=0}^n \scf(G^{m-1})\otimes \scf(G^{n-m-1})\\
\psi & \mapsto & \dd\sum_{A\subseteq \{1,2,\ldots,n\}}\Down_{(A,\overline{A})}^{(\{1,2,\ldots,n\})} [\tau,\alpha,\beta](\psi),
\end{array}\end{equation}
linearly, where by convention if $A$ or $\overline{A}$ is empty, then we get  $\chi^{\emptyset}\otimes \psi$ and  $\psi\otimes \chi^{\emptyset}$, respectively.  By Lemma \ref{LemmaCoAssociativity} (b) this product is co-associative.

If $\psi=\psi_1\otimes \psi_2\otimes \cdots \otimes \psi_n\in \scf(G^{n-1})$ factors, then for $1\leq j\leq n-1$,
\begin{equation}\label{DownFormulas}
\Down_{(A,\overline{A})}^{(\{1,2,\ldots,n\})} [\tau,\alpha,\beta](\psi)_j=\left\{\begin{array}{ll}
\psi_j & \text{if $\bc(A,\overline{A})_j=1$,}\\
\langle\psi_j,\alpha\rangle \tau & \text{if $\bc(A,\overline{A})_j=0$, $\llc(A,\overline{A})_j=\lc(A,\overline{A})_j=1$,}\\
\langle\psi_j,\alpha\rangle \chi^{()} & \text{if $\bc(A,\overline{A})_j=\lc(A,\overline{A})_j=0$, $\llc(A,\overline{A})_j=1$,}\\
\langle\psi_j,\beta\rangle \tau & \text{if $\bc(A,\overline{A})_j=\llc(A,\overline{A})_j=0$, $\lc(A,\overline{A})_j=1$,}\\
\langle\psi_j,\beta\rangle \chi^{()} & \text{if $\bc(A,\overline{A})_j=\lc(A,\overline{A})_j=\llc(A,\overline{A})_j=0$.}\\
\end{array}\right.
\end{equation}
Note that this expresses the product without reshuffling into two factors, so for example, if $A=\{2,3,6,7,8\}\subseteq \{1,2,\ldots, 10\}$, then
\begin{align*}
&\Down_{(A,\overline{A})}^{(\{1,2,\ldots,10\})} [\tau,\alpha,\beta](\psi)\\
&= \tau_{(A,\overline{A})}\odot \Inf_{G^{9}_{\bc(A,\overline{A})}}^{G^{9}_{\lc(A,\overline{A})}}\circ \Def_{G^{9}_{\bc(A,\overline{A})}}^{G^{9}}\Big(\underset{{\color{blue}1}}{ {\color{white}\otimes}}\beta\odot\psi_1{\color{red}\overset{2}{\otimes}}\psi_2{\color{red}\overset{3}{\otimes}}\alpha\odot\psi_3{\color{blue}\underset{4}{\otimes}}\psi_4{\color{blue}\underset{5}{\otimes}}\beta\odot\psi_5{\color{red}\overset{6}{\otimes}}\psi_6{\color{red}\overset{7}{\otimes}}\psi_7{\color{red}\overset{8}{\otimes}}\alpha\odot\psi_8{\color{blue}\underset{9}{\otimes}}\psi_9{\color{blue}\underset{10}{\otimes}}\chi^{()}\Big)\\
&=\tau_{(A,\overline{A})}\odot \Inf_{G^{9}_{\bc(A,\overline{A})}}^{G^{9}_{\lc(A,\overline{A})}}\Big(\underset{{\color{blue}1}}{ {\color{white}\otimes}}\langle\psi_1,\beta\rangle\chi^{()}{\color{red}\overset{2}{\otimes}}\psi_2{\color{red}\overset{3}{\otimes}}\langle\psi_3,\alpha\rangle\chi^{()}{\color{blue}\underset{4}{\otimes}}\psi_4{\color{blue}\underset{5}{\otimes}}\langle\psi_5,\beta\rangle\chi^{()}{\color{red}\overset{6}{\otimes}}\psi_6{\color{red}\overset{7}{\otimes}}\psi_7{\color{red}\overset{8}{\otimes}}\langle\psi_8,\alpha\rangle\chi^{()}{\color{blue}\underset{9}{\otimes}}\psi_9{\color{blue}\underset{10}{\otimes}}\chi^{()}\Big)\\
&=\underset{{\color{blue}1}}{ {\color{white}\otimes}}\langle\psi_1,\beta\rangle\tau{\color{red}\overset{2}{\otimes}}\psi_2{\color{red}\overset{3}{\otimes}}\langle\psi_3,\alpha\rangle\tau{\color{blue}\underset{4}{\otimes}}\psi_4{\color{blue}\underset{5}{\otimes}}\langle\psi_5,\beta\rangle\tau{\color{red}\overset{6}{\otimes}}\psi_6{\color{red}\overset{7}{\otimes}}\psi_7{\color{red}\overset{8}{\otimes}}\langle\psi_8,\alpha\rangle\chi^{()}{\color{blue}\underset{9}{\otimes}}\psi_9{\color{blue}\underset{10}{\otimes}}\chi^{()}
\end{align*}
which unshuffles to
\[\Big({\color{red}\overset{2}{{\color{white}\otimes}}}\psi_2{\color{red}\overset{3}{\otimes}}\langle\psi_3,\alpha\rangle\tau{\color{red}\overset{6}{\otimes}}\psi_6{\color{red}\overset{7}{\otimes}}\psi_7{\color{red}\overset{8}{\otimes}}\langle\psi_8,\alpha\rangle\chi^{()}\Big)\otimes\Big(
\underset{{\color{blue}1}}{ {\color{white}\otimes}}\langle\psi_1,\beta\rangle\tau{\color{blue}\underset{4}{\otimes}}\psi_4{\color{blue}\underset{5}{\otimes}}\langle\psi_5,\beta\rangle\tau{\color{blue}\underset{9}{\otimes}}\psi_9{\color{blue}\underset{10}{\otimes}}\chi^{()}\Big)\in {\color{red}\scf(G^{5-1})}\otimes {\color{blue}\scf(G^{5-1})}.
\]
In this case, if we are working in a specific basis we need to know how to decompose, $\alpha$, $\beta$ and $\tau$.

\begin{runningexample}
In our example with $\scf(G)=\CC\spanning\{\mathds{1},\reg\}$, consider two cases.  
\begin{description}
\item[Case $\tau=\alpha=\beta=\mathds{1}$.]  For $A\subseteq \{1,2,\ldots,n\}$, $\underline{b}\in \{0,1\}^{n-1}$,
\begin{equation*}
\Down_{(A,\overline{A})}^{(\{1,\ldots,n\})}[\mathds{1},\mathds{1},\mathds{1}](\chi^{\underline{b}})=\left\{\begin{array}{ll}
0 & \text{if $\bc(A,\overline{A})_j=0=b_j$ for some $1\leq j\leq n-1$},\\
\chi^{\underline{b}_A}\otimes \chi^{\underline{b}_{\overline{A}}} & \text{if $\bc(A,\overline{A})_j=0$ implies $b_j=1$,} 
\end{array}\right.
\end{equation*}
where if $A=\{a_1<a_2<\ldots <a_\ell\}$, then $\underline{b}_A\in \{0,1\}^{\ell-1}$ is the sequence given by 
\[
(\underline{b}_A)_j=b_{a_j}.
\]
Note that if we encode $\{0,1\}^{n-1}$ as integer compositions of $n$ by letting $1$'s denote carriage returns, then this operation deshuffles the parts.

\item[Case $\tau=\reg$, $\alpha=\mathds{1}$, $\beta=(\reg-\mathds{1})^*$.] For $A\subseteq \{1,2,\ldots,n\}$, $\underline{b}\in \{0,1\}^{n-1}$,
\begin{align*}
\Down_{(A,\overline{A})}^{(\{1,\ldots,n\})}[\reg,&\mathds{1},(\reg-\mathds{1})^*](\chi^{\underline{b}})\\
&=\left\{\begin{array}{ll}
0 & \text{if $\llc(A,\overline{A})_j<b_j$ for some $j$},\\
\dd\sum_{{b'\in \{0,1\}^A,b''\in \{0,1\}^{\overline{A}}\atop b'_a=b_a\text{ if $a\in A$, $\bc(A,\overline{A})_a=1$}}\atop b''_a=b_a\text{ if  $a\in \overline{A}$, $\bc(A,\overline{A})_a=1$}}\chi^{b'}\otimes \chi^{b''} & \text{if $\llc(A,\overline{A})_j=0$ implies $b_j=0$.} 
\end{array}\right.
\end{align*}
\end{description}
\end{runningexample}

\subsection{Hopf compatibility}

As defined in (\ref{HopfProduct}) and (\ref{HopfCoProduct}), the product and coproduct generally do not give a bi-algebra, but under some minor restrictions on $\iota$, $\tau$, $\alpha$, and $\beta$  we in fact not only get a Hopf algebra, but the underlying structure for a Hopf monoid in species.  Our main result in this section is as follows, where for a set of integers $X$ we let
\[
X_{\leq m}=\{x\in X\mid x\leq m\}\quad \text{and}\quad X_{>m}=\{x\in X\mid x>m\}.
\]
\begin{theorem}\label{HopfCompatibility}
Let $\iota,\tau, \alpha,\beta\in \scf(G)$.  For all $(A,B)\vDash \{1,2,\ldots, n\}$,  $\psi\in \scf(G^{|A|-1})$, and $\gamma\in \scf(G^{|B|-1})$,
\begin{align*}
&\Down_{(A,B)}^{(\{1,2,\ldots,n\})}[\tau,\alpha,\beta]\circ\Inf_{(\{1,2,\ldots,|A|\},\{|A|+1,\ldots, n\})}^{(\{1,2,\ldots, n\})}[\iota](\psi\otimes \gamma)\\
&=\Big(\Inf_{(A_{\leq |A|},A_{>|A|})}^{(A)}[\iota]\otimes \Inf_{(B_{\leq |A|},B_{>|A|})}^{(B)}[\iota]
\Big)
\Big(\Down_{(A_{\leq |A|},B_{\leq |A|})}^{(\{1,\ldots,|A|\})}[\tau,\alpha,\beta](\psi)\otimes \Down_{(A_{> |A|},B_{> |A|})}^{(\{|A|+1,\ldots,n\})}[\tau,\alpha,\beta](\gamma)\Big)
\end{align*}
if and only if $\iota=\tau$ and  $\langle\iota,\alpha\rangle=\langle\iota,\beta\rangle=1$.
\end{theorem}

\begin{proof}
Since we have a basis whose elements factor, we may assume WLOG that $\psi$ and $\gamma$ factor.  By Lemmas \ref{LemmaAssociativity} (a) and \ref{LemmaCoAssociativity} (a) we may evaluate this equality coordinate-wise.  First, consider the case $j=m=|A|$.  Then	
\begin{align*}
\Down_{(A,B)}^{(\{1,2,\ldots,n\})}[\tau,\alpha,\beta]&\circ\Inf_{(\{1,2,\ldots,m\},\{m+1,\ldots, n\})}^{(\{1,2,\ldots, n\})}[\iota](\psi\otimes \gamma)_{m} \\
&= \left\{\begin{array}{ll}\iota  & \text{if $\bc(A,B)_{m}=1$,}\\
\langle \iota,\beta\rangle \tau & \text{if $\bc(A,B)_{m}=0$, $\llc(A,B)_{m}=0$, $\lc(A,B)_{m}=1$,}\\
\langle \iota,\alpha\rangle \tau & \text{if $\bc(A,B)_{m}=0$, $\llc(A,B)_{m}=1$, $\lc(A,B)_{m}=1$,}\\
\langle \iota,\beta\rangle  \chi^{()}& \text{if $\bc(A,B)_{m}=0$, $\llc(A,B)_{m}=0$, $\lc(A,B)_{m}=0$,}\\
\langle \iota,\alpha\rangle \chi^{()} & \text{if $\bc(A,B)_{m}=0$, $\llc(A,B)_{m}=1$, $\lc(A,B)_{m}=0$.}\\
\end{array}
\right.
\end{align*}
On the other hand, 
\begin{align*}
\Big(\Inf_{(A_{\leq m},A_{>m})}^{(A)}[\iota]\otimes& \Inf_{(B_{\leq m},B_{>m})}^{(B)}[\iota]
\Big)
\Big(\Down_{(A_{\leq m},B_{\leq m})}^{(\{1,\ldots,m\})}[\tau,\alpha,\beta](\psi)\otimes \Down_{(A_{> m},B_{> m})}^{(\{m+1,\ldots,n\})}[\tau,\alpha,\beta](\gamma)\Big)\\
&= \left\{\begin{array}{ll}\iota  & \text{if $\lc(A,B)_{m}=1$,}\\
\chi^{()} & \text{if $\lc(A,B)_{m}=0$.}
\end{array}
\right.
\end{align*}
If we choose $A$ so that $\lc(A,B)_{m}=0$ we get equality only if
\[
\langle \iota,\beta\rangle=1.  
\]
Next choose $A$ such that $\lc(A,B)_{m}=1$, $\bc(A,B)_{m}=0$ and $\llc(A,B)_{m}=0$.  Then we get equality only if
\[
\iota=\langle \iota, \beta\rangle \tau=\tau.
\]
If $A$ satisfies  $\lc(A,B)_{m}=1$, $\bc(A,B)_{m}=0$ and $\llc(A,B)_{m}=1$, we get equality only if
\[
\iota=\langle \iota, \alpha\rangle \tau= \langle \iota, \alpha\rangle\iota,
\]
so $\langle\iota,\alpha\rangle=1$.   In particular, we have established that we get equality only if $\iota=\tau$ and  $\langle\iota,\alpha\rangle=\langle\iota,\beta\rangle=1$.  Furthermore, if these conditions hold we get equality when $j=m$.

Suppose $\iota=\tau$,  $\langle\iota,\alpha\rangle=\langle\iota,\beta\rangle=1$, and $j\neq m$. Then
\begin{align*}
\Down_{(A,B)}^{(\{1,2,\ldots,n\})}[\tau,\alpha,\beta]&\circ\Inf_{(\{1,2,\ldots,m\},\{m+1,\ldots, n\})}^{(\{1,2,\ldots, n\})}[\iota](\psi\otimes \gamma)_{j} \\
&= \left\{\begin{array}{ll} (\psi\otimes \gamma)_j  & \text{if $\bc(A,B)_{j}=1$,}\\
\langle (\psi\otimes \gamma)_j,\beta\rangle \iota & \text{if $\bc(A,B)_{j}=0$, $\llc(A,B)_{j}=0$, $\lc(A,B)_{j}=1$,}\\
\langle (\psi\otimes \gamma)_j,\alpha\rangle \iota & \text{if $\bc(A,B)_{j}=0$, $\llc(A,B)_{j}=1$, $\lc(A,B)_{j}=1$,}\\
\langle (\psi\otimes \gamma)_j,\beta\rangle  \chi^{()}& \text{if $\llc(A,B)_{j}=0$, $\lc(A,B)_{j}=0$,}\\
\langle (\psi\otimes \gamma)_j,\alpha\rangle  \chi^{()}& \text{if $\llc(A,B)_{j}=1$, $\lc(A,B)_{j}=0$.}
\end{array}
\right.
\end{align*}
On the other hand, if $j\in X$ with $X\in \{(A_{\leq m},B_{\leq m}),(A_{>m},B_{>m})\}$, then
\begin{align*}
\Big(\Inf_{(A_{\leq m},A_{>m})}^{(A)}[\iota]\otimes& \Inf_{(B_{\leq m},B_{>m})}^{(B)}[\iota]
\Big)
\Big(\Down_{(A_{\leq m},B_{\leq m})}^{(\{1,\ldots,m\})}[\tau,\alpha,\beta](\psi)\otimes \Down_{(A_{> m},B_{> m})}^{(\{m+1,\ldots,n\})}[\tau,\alpha,\beta](\gamma)\Big)\\
&= \left\{\begin{array}{ll} (\psi\otimes \gamma)_j  & \text{if $\bc(X)_{j}=1$,}\\
\langle(\psi\otimes \gamma)_j,\beta\rangle \iota  & \text{if $\bc(X)_j=0$,  $\llc(X)_j=0$, $\lc(X)_{j}=1$, }\\
\langle(\psi\otimes \gamma)_j,\alpha\rangle \iota  & \text{if $\bc(X)_j=0$,   $\llc(X)_j=1$, $\lc(X)_{j}=1$,}\\
\langle (\psi\otimes \gamma)_j,\beta\rangle \iota & \text{if $\lc(X)_{j}=0$, $\llc(X)_j=0$, $\lc(A,B)_j=1$,}\\
\langle (\psi\otimes \gamma)_j,\alpha\rangle \iota & \text{if $\lc(X)_{j}=0$, $\llc(X)_j=1$, $\lc(A,B)_j=1$,}\\
\langle (\psi\otimes \gamma)_j,\beta\rangle\chi^{()} & \text{if $\llc(X)_{j}=0$, $\lc(A,B)_j=0$,}\\
\langle (\psi\otimes \gamma)_j,\alpha\rangle\chi^{()} & \text{if $\llc(X)_{j}=1$, $\lc(A,B)_j=0$.}
\end{array}
\right.
\end{align*}
Since for $j\neq m$, $\bc(X)_j=\bc(A,B)_j$ and if $\lc(X)_j=1$, then $\llc(X)_j=\llc(A,B)_j$, we can check that we get equality.
\end{proof}

\begin{remark}
Theorem \ref{HopfCompatibility} implies that we not only get a Hopf algebra, but that there is in fact an underlying Hopf monoid.
\end{remark}

We say $(\iota,\alpha,\beta)$ is a \emph{Hopf triple} if $\iota,\alpha,\beta\in \scf(G)$ with $\langle\iota,\alpha\rangle=\langle\iota,\beta\rangle=1$.  For a Hopf triple, let $\cH_{(\iota,\alpha,\beta)}$ denote the Hopf algebra
\[
\cH_{(\iota,\alpha,\beta)}=\bigoplus_{n\geq 0} \scf(G^{n-1})
\]
whose product (\ref{HopfProduct}) uses $\iota=\iota$ and coproduct (\ref{HopfCoProduct}) uses $\tau=\iota$, $\alpha=\alpha$ and $\beta=\beta$.

\subsection{Freeness}

All the above Hopf algebras are free, as proved in the next result.

\begin{proposition}\label{HopfIsFree}
Let $(\iota,\alpha,\beta)$ be a Hopf triple. Then the Hopf algebra $\cH_{(\iota,\alpha,\beta)}$ is free.
\end{proposition}
\begin{proof}
Fix a basis $\cB\subseteq \scf(G)$ with $\iota\in \cB$, and let $\tilde\cB=\cB-\{\iota\}$.  For $n\in\ZZ_{\geq 1}$,  the set $\cB^{\otimes (n-1)}$ is a basis for $\scf(G^{n-1})$, and let   
\[
\tilde\cB^\all=\bigcup_{n\geq 0} \tilde\cB^{\otimes (n-1)},
\quad\text{where}\quad 
\tilde\cB^{\otimes (n-1)}=\{\psi_1\otimes\cdots\otimes \psi_{n-1}\otimes \chi^{()}\mid \psi_1,\ldots,\psi_{n-1}\in \tilde\cB\}.
\]
If $\psi_1\otimes \cdots \otimes \psi_{n-1}\otimes \chi^{()}\in \cB^{\otimes (n-1)}$, let $(i_1,i_2,\ldots, i_\ell)$ be the subsequence of $(1,2,\ldots,n-1)$ recording the indices with $\psi_{i_j}=\iota$.   Then
\begin{align*}
\psi_1\otimes &\cdots \otimes \psi_{n-1}\otimes \chi^{()}\\
 &=
\Big(\psi_1\otimes\cdots \otimes \psi_{i_1-1}\otimes \chi^{()}\Big)\Big(\psi_{i_1+1}\otimes \cdots\otimes \psi_{i_2-1}\otimes \chi^{()}\Big)\cdots \Big(\psi_{i_{\ell-1}+1}\otimes \cdots\otimes \psi_{n-1}\otimes \chi^{()}\Big)
\end{align*}
is a unique factorization into a product of elements in $\tilde\cB^\all$.  Thus, $\tilde\cB^\all$ is a set of free generators for $\cH_{(\iota,\alpha,\beta)}$.
\end{proof}

\begin{runningexample}
If $\scf(G)=\CC\spanning\{\mathds{1},\reg\}$ and $\iota\in \{\mathds{1},\reg\}$, then $|\tilde\cB|=1$.  Thus, each degree has exactly one generator.
\end{runningexample}

\section{Favorite Examples: $\NSym$ and $\QSym$}\label{SectionExamples}

The main motivation for this work was understanding the different representation theoretic interpretations for $\NSym$ the authors had come across in recent years.  This section focuses  on the case where $\scf(G)=\CC\spanning\{\mathds{1},\reg\}$, so that $\dim(\scf(G^{n-1}))=2^{n-1}$.

\subsection{Free cocommutative graded Hopf algebras}

Let $\cH=\bigoplus_{n\geq 0}\cH_n$ be a free graded cocommutative Hopf algebra.  By Milnor--Moore, it is generated by its set of primitive elements $\prim(\cH)$.  
\begin{lemma}
Let $\cH=\bigoplus_{n\geq 0}\cH_n$ be a free graded cocommutative Hopf algebra.  Then $\cH$ is freely generated by a subset of $\prim(\cH)$.
\end{lemma}
\begin{proof}
Let $F$ be a set that freely generates $H$ such that each generator has homogeneous degree, and let $F_{\leq n}= F\cap \bigoplus_{j\leq n} \cH_j$.  We prove by inducting on $n$ that $F_{\leq n}$ may be replaced by $F_{\leq n}'\subseteq \prim(\cH_n)$ such that the set $F'$ obtained by replacing $F_{\leq n}$ in $F$ by $F_{\leq n}'$ still freely generates $\cH$. Note that for $n=1$, $F_{\leq 1}\subseteq\prim(\cH)$.  

Suppose we can replace $F_{\leq (n-1)}$ by $F'_{\leq (n-1)}\subseteq \prim(\cH)$.   If $F_n\subseteq \prim(\cH)$, then we are done.  Else, let $f\in F_n-\prim(\cH)$.  Since
$\prim(\cH)$ generates $\cH$, 
\[f=p_{<n}+p_{=n},\]
where $p_{<n}$ is a linear combination of products of primitives from a lower degree (WLOG in $F'_{\leq (n-1)}$), and $p_{=n}$ is a linear combination of primitives of degree $n$.  But
\[
f-p_{<n}\in \prim(\cH),
\]
so if $F''$ is the set where we replace $f$ in $F_n$ by $f-p_{<n}$, then $F''$ also generates $\cH$.   To see that $F''$ is still free suppose there is a polynomial $Y$ in $F''$ involving $f-p_{<n}$ that is zero, so
\[
Y(f-p_{<n})=0.
\]
However, $f-p_{<n}$ is a polynomial in the elements of $F$, so by the freeness of $F$, $Y$ must be identically zero.  Thus $F''$ is still free.  In this way, every element in $F_n-\prim(\cH)$ can be replaced by a primitive, and we obtain a new set $F_n'$ of primitives such that $F'_{\leq n}=F_{\leq (n-1)}'\cup F_n'$ gives us a desired replacement.
\end{proof}

We deduce the main result of this section. 

\begin{theorem}
Let $\cH=\bigoplus_{n\geq 0}\cH_n$ and $\cH'=\bigoplus_{n\geq 0}\cH'_n$ be two free graded cocommutative Hopf algebra. Then $\cH\cong\cH'$ if and only if
\[\dim(\cH_n)=\dim(\cH'_n)\quad \text{for all $n\geq 0$}.\]
\end{theorem}
\begin{proof}
Let $F\subseteq \prim(\cH)$ and $F'\subseteq \prim(\cH')$ be sets of homogeneous free generators.  Suppose $n$ is minimal such that $|F_n|\neq |F_n'|$.  Then
\begin{align*}
\#\{&x_1\cdots x_\ell\mid x_1,\ldots, x_\ell\in \bigcup_{j<n} F_j, \deg(x_1)+\cdots+\deg(x_\ell)=n,\ell\geq 1\}+|F_n|=\dim(\cH_n)\\
&=\dim(\cH_n')=\#\{x_1\cdots x_\ell\mid x_1,\ldots, x_\ell\in \bigcup_{j<n} F'_j, \deg(x_1)+\cdots+\deg(x_\ell)=n,\ell\geq 1\}+|F'_n|.
\end{align*}
Since $|F_j|=|F_j'|$ for all $j<n$, we conclude that $|F_n|=|F_n'|$, contradicting our choice of $n$.  Thus, $|F_n|=|F_n'|$ for all $n$, and we may fix a graded bijection $\phi:F\rightarrow F'$.  Then $\phi$ extends multiplicatively to an algebra homomorphism $\Phi$.  Suppose $f_1,\ldots, f_\ell\in F$.  Then
\begin{align*}
\Delta(\Phi(f_1\cdots f_\ell))&=\Delta(\Phi(f_1))\cdots \Delta(\Phi(f_\ell))\\
&=\sum_{A\sqcup B=\{1,2,\ldots, \ell\}} \Phi(f_A)\otimes \Phi(f_B)\\
&=\Phi(\Delta(f_1\cdots f_\ell)),
\end{align*}
where $f_A=f_{a_1}\cdots f_{a_{|A|}}$ if $A=\{a_1<a_2<\cdots<a_{|A|}\}$.  Extend linearly to get that $\Phi$ is also a co-algebra isomorphism.
\end{proof}

\begin{corollary} \label{ItsAllNSym}  
Fix a supercharacter theory $(\Ch,\Cl)$ of $G$ and let $\iota,\alpha,\beta\in \scf(G)$ form a Hopf triple.
\begin{enumerate}
\item[(a)]  If $(\tilde\Ch,\tilde\Cl)$ is another supercharacter theory of $G$ with $\dim(\scf(G))=\dim(\widetilde{\scf}(G))$, and  $\tilde\iota,\tilde\alpha,\tilde\beta\in \widetilde\scf(G)$ forms a Hopf triple, then
\[
\cH_{(\iota,\alpha,\beta)}\cong \cH_{(\tilde\iota,\tilde\alpha,\tilde\beta)}.
\]
\item[(b)] If $\dim(\scf(G))=2$, then 
\[\cH_{(\iota,\alpha,\beta)}\cong \NSym.\]
\end{enumerate}
\end{corollary}

\subsection{$\NSym$} \label{NSym}

While we obtain a number of different versions of $\NSym$ via Corollary \ref{ItsAllNSym} the choices naturally pick out certain bases.    For $\tau,\iota\in \scf(G)$ and $\mu\vDash n$, let
$(\tau)^\iota_\mu\in\scf(G^{n-1})$ be given by
\[
(\tau)^\iota_\mu=\tau^{\otimes (\mu_1-1)}\otimes \iota\otimes\tau^{\otimes (\mu_2-1)}\otimes\iota\otimes \cdots \otimes \iota\otimes \tau^{\otimes(\mu_\ell-1)}\otimes \chi^{()} \in\scf(G^{n-1}).
\]
In fact, since
\[(\tau)^\iota_{(n)}= \tau^{\otimes(n-1)}\otimes \chi^{()} \]
in $\cH_{(\iota,\alpha,\beta)}$ we have
\begin{equation}\label{SimpleFactoring}
(\tau)_\mu^\iota =(\tau)_{(\mu_1)}^\iota(\tau)_{(\mu_2)}^\iota\cdots (\tau)_{(\mu_\ell)}^\iota.
\end{equation}
By construction, we have
\[(\tau)^\iota_\mu=(\iota)^\tau_{\mu'},\]
where $\mu'$ is the conjugate ribbon where we replace right steps with down steps and vice-versa.  In terms of binary sequences, conjugation corresponds to switching zeroes and ones.  

In particular, if $\scf(G)=\CC\spanning\{\mathds{1},\reg\}$, then the supercharacters are given by
\[
\{(\reg-\mathds{1})^{\mathds{1}}_\mu\mid\mu\vDash n\}=\{(\mathds{1})^{\reg-\mathds{1}}_\mu\mid\mu\vDash n\},
\]
and the superclass identifiers are given by
\[
\Big\{\Big(\frac{\reg}{|G|}\Big)^{\mathds{1}}_\mu\mid \mu\vDash n\Big\}=\{(\mathds{1})^{\reg/|G|}_\mu\mid \mu\vDash n\}.
\]

Note that if $(\iota,\alpha,\beta)$ is a Hopf triple and $\alpha\in \CC\spanning\{\beta\}$, then due to the compatibility conditions with $\iota$, we have $\alpha=\beta$.  If $\alpha\neq \beta$, let $\{\alpha^\dual,\beta^\dual\}$ denote the dual basis to $\{\alpha,\beta\}$ with respect to the usual inner product $\langle\cdot,\cdot\rangle$, so that $\iota=\alpha^\dual+\beta^\dual$.

\begin{lemma}\label{ExplicitNSymIsos}
Suppose $\scf(G)=\CC\spanning\{\mathds{1},\reg\}$ and $(\iota,\alpha,\beta)$ is a Hopf triple.
\begin{description}
\item[Case $\alpha=\beta$.] Let $\tau\in \scf(G)$ be a nonzero element such that $\langle\tau,\beta\rangle=0$ and $\CC\spanning\{\tau,\iota\}=\scf(G)$.  Then $\{(\tau)^\iota_{(n)}\mid n\in \ZZ_{\geq 1}\}$ is a set of primitives that freely generate $\cH_{(\iota,\alpha,\beta)}$ and
\[
\begin{array}{ccc}
\cH_{(\iota,\alpha,\beta)} & \longrightarrow & \NSym\\
(\tau)_\mu^\iota & \mapsto & \langle\cdot,p_\mu\rangle,
\end{array}
\]
is a Hopf algebra isomorphism, where $p_\mu$ is a shuffle basis element of $\QSym$.
\item[Case $\alpha\neq \beta$.]   The set $\{(\alpha^\dual)^\iota_{(n)}\mid n\in \ZZ_{\geq 1}\}$ freely generates $\cH_{(\iota,\alpha,\beta)}$ and
\[
\begin{array}{ccc}
\cH_{(\iota,\alpha,\beta)} & \longrightarrow & \NSym\\
(\alpha^\dual)_\mu^\iota & \mapsto & H_\mu,
\end{array}
\]
is a Hopf algebra isomorphism.  In this case, 
\[
(\alpha^\dual)^{\beta^\dual}_\mu\mapsto R_\mu.
\]
\end{description}
\end{lemma}

\begin{proof}
In both cases,  (\ref{SimpleFactoring}) implies that the sets generate freely, using a similar argument as in the proof of Proposition \ref{HopfIsFree}. For the coproduct
\begin{align*}
\Delta((\tau)_{(n)}^\iota )&= (\tau)_{(n)}^\iota\otimes \chi^\emptyset+  \chi^\emptyset\otimes (\tau)_{(n)}^\iota, \\
\Delta((\alpha^\dual)^{\iota}_{(n)})&=\sum_{j=0}^n (\alpha^\dual)_{(j)}^\iota \otimes (\alpha^\dual)_{(n-j)}^\iota, 
\end{align*}
and the results follow by the multiplicativity of the coproduct.
\end{proof}

\begin{examples}
If $\iota=\alpha=\beta=\mathds{1}$, then by letting $\tau=\reg-\mathds{1}$, we get the that the supercharacters 
\[
(\tau)^\iota_\mu=\chi^{\binaryzero(\mu)}
\]
 are dual to the shuffle basis of $\QSym$.

 Let $\iota=\reg$, $\alpha=\mathds{1}$ and $\beta=\frac{\reg-\mathds{1}}{|G|-1}$. Then $\mathds{1}^\dual=\mathds{1}$ and the permutation characters 
\[(\mathds{1})^\reg_\mu=\Ind_{G^{n-1}_{\binaryone(\mu)}}^{G^{n-1}}(\mathds{1})\] 
correspond to the $H$-basis.
In this case, the supercharacters correspond to the ribbon basis of $\NSym$.
\end{examples}

\subsection{The relation of $\NSym$ to other Hopf algebras}\label{Embeddings}

Recall, that $\NSym$ injects into the Hopf algebra $\NCSym$ of symmetric functions in noncommuting variables.  There are various ways to realize this inclusion, but the following gives one.

\begin{proposition}
Let $G=\FF_q^+$.  Then the function $\Inf:  \cH_{(\mathds{1},\reg,\reg)}  \longrightarrow  \scf(\UT_\all)$ given by extending
\[
\Inf_{\UT_n/[\UT_n,\UT_n]}^{\UT_n}(\psi) \quad \text{for $\psi\in \scf(G^{n-1})$,}
\]
linearly, is an injective Hopf algebra morphism that sends the supercharacters of $\scf(G^\all)$ to the supercharacters containing the linear characters of $\scf(\UT_\all)$.
\end{proposition}

There is also (at least) one injective homomorphism into the Malvenuto--Reutenauer Hopf algebra.  The paper \cite{ATmr} relates the Malvenuto--Reutenauer algebra to the representation theory of 
\[
\fkut_n=\{u\in M_n(\FF_q)\mid u_{ji}=0, \text{ if } i\leq j\}.
\]
For $w\in S_n$, define 
\begin{align*}
\fkut_w&=\Big\{u\in \fkut_n\mid u_{ij}=0, \text{ if } j-i> \#\{h<i\mid w(h)>w(i)\}\Big\} \\
\fkut_w^+&=\Big\{u\in \fkut_n\mid u_{ij}=0, \text{ if } j-i-1> \#\{h<i\mid w(h)>w(i)\}\Big\}. 
\end{align*}
For example, if $w=(4,2,9,5,8,3,1,7,6)$,
\[
\fkut_w=\left[\begin{array}{ccccccccc}
0 & \ast & \ast & \ast & \ast & \ast & \ast & 0 & 0 \\ 
\cdot  & 0  & \ast & 0 & 0 & 0 & 0 & 0 & 0 \\ 
\cdot  & \cdot  & 0  & \ast & \ast & \ast & \ast & 0 & 0  \\ 
\cdot  & \cdot  & \cdot  & 0  & 0 & 0 & 0 & 0 & 0  \\ 
\cdot  & \cdot  & \cdot  & \cdot  & 0  & \ast & 0 & 0 & 0   \\ 
\cdot  & \cdot  & \cdot  & \cdot  & \cdot  & 0  & \ast & \ast & \ast    \\ 
\cdot  & \cdot  & \cdot  & \cdot  & \cdot  & \cdot  & 0  & \ast & \ast     \\
\cdot  & \cdot  & \cdot  & \cdot  & \cdot  & \cdot  & \cdot  & 0  & \ast    \\ 
\cdot  & \cdot  & \cdot  & \cdot  & \cdot  & \cdot  & \cdot  & \cdot  & 0    
\end{array}\right]\quad \text{and}\quad \fkut_w ^+=\left[\begin{array}{ccccccccc}
0 & \ast & \ast & \ast & \ast & \ast & \ast & \color{red}\ast & 0 \\ 
\cdot  & 0  & \ast &\color{red} \ast & 0 & 0 & 0 & 0 & 0 \\ 
\cdot  & \cdot  & 0  & \ast & \ast & \ast & \ast & \color{red}\ast & 0  \\ 
\cdot  & \cdot  & \cdot  & 0 & \color{red}\ast & 0 & 0 & 0 & 0  \\ 
\cdot  & \cdot  & \cdot  & \cdot  & 0  & \ast & \color{red}\ast & 0 & 0   \\ 
\cdot  & \cdot  & \cdot  & \cdot  & \cdot  & 0  & \ast & \ast & \ast    \\ 
\cdot  & \cdot  & \cdot  & \cdot  & \cdot  & \cdot  & 0  & \ast & \ast     \\
\cdot  & \cdot  & \cdot  & \cdot  & \cdot  & \cdot  & \cdot  & 0  & \ast    \\ 
\cdot  & \cdot  & \cdot  & \cdot  & \cdot  & \cdot  & \cdot  & \cdot  & 0    
\end{array}\right].
\]

For $w\in S_n$, let $\tilde{w}$ be the set composition of $\{1,2,\ldots, n\}$ given by
\[
\tilde{w}=(\{w(1)\},\{w(2)\},\ldots, \{w(n)\}).
\]
Then we have the following composition of functions
\[
\begin{array}{r@{\ }c@{\ }c@{\ }c@{\ }c@{\ }c@{\ }c@{\ }c}
\twst_w: & \scf(G^{n-1}) & \longrightarrow & \scf(\fkut_w) & \longrightarrow & \scf(\fkut_w^+) & \longrightarrow & \scf(\fkut_n)\\
 & \psi & \mapsto &\Inf_{\{1\}}^{\fkut_w}\circ \Down_{\widetilde{w}} [\reg,\mathds{1},(\reg-\mathds{1})^\dual](\psi)\\
 & & & \gamma & \mapsto & \Inf[\reg-\mathds{1}](\gamma)\\
 & & & & &  \theta & \mapsto & \Ind_{\fkut_w}^{\fkut_n}(\theta).
\end{array}
\]

If we apply the first function to a supercharacter $(\mathds{1})^{\reg-\mathds{1}}_\mu$, we obtain
\[
\Inf_{\{1\}}^{\fkut_w}\circ \Down_{\widetilde{w}} [\mathds{1},(\reg-\mathds{1})^\dual]\Big((\mathds{1})^{\reg-\mathds{1}}_\mu\Big)=0,
\]
unless $w^{-1}(i)>w^{-1}(i+1)$ if and only if $i\in \{ \mu_1,\mu_1+\mu_2,\ldots, n-\mu_\ell\}.$  Thus,
\[
\twst_w\Big((\mathds{1})^{\reg-\mathds{1}}_\mu\Big)=\left\{\begin{array}{ll}
\chi^{w} & \text{if $w^{-1}(i)>w^{-1}(i+1)$ if and only if $i\in \{ \mu_1,\mu_1+\mu_2,\ldots, n-\mu_\ell\}.$,}\\ 0 & \text{otherwise}.
\end{array}\right.
\]
Extend linearly, to get the function
\[
\begin{array}{rccc}
\twst: & \scf(G^{n-1}) & \longrightarrow & \scf(\fkut_n)\\
& \psi & \mapsto & \dd\sum_{w\in S_n} \twst_w(\psi).
\end{array}
\]

\begin{proposition}
The function $\twst:\cH_{(\reg,\mathds{1},(\reg-\mathds{1})^\dual)}\longrightarrow\scf(\fkut_\all)$ is an injective Hopf algebra morphism.
\end{proposition}
\begin{proof}
First note that for $\mu\vDash n$,
\[
\twst\Big((\mathds{1})^{\reg-\mathds{1}}_\mu \Big)=\sum_{ w\in  S_n\atop \desc(w^{-1})=\{\mu_1,\ldots, n-\mu_\ell\}} \chi^w\leftrightarrow R_\mu,
\]
where for $v\in S_n$,
\[
\desc(v)=\{1\leq i\leq n\mid v(i)>v(i+1)\}.
\]

Thus, the composition
\[
\begin{array}{ccccccccc}
\NSym & \longrightarrow &\cH_{(\reg,(\reg-\mathds{1})^\dual,\mathds{1})} & \longrightarrow & \scf(\fkut_\all) & \longrightarrow & \MR\\
R_{\mu} & \mapsto & (\mathds{1})^{\reg-\mathds{1}}_\mu& \mapsto & \twst \Big( (\mathds{1})^{\reg-\mathds{1}}_\mu  \Big) & \mapsto & R_{\mu}
\end{array}
\]
gives the usual embedding of $\NSym$ into $\MR$ \cite{GR}. 
\end{proof}

\section{Combinatorial Hopf algebras} \label{CombinatorialHopfAlgebras}

We return to the general setting, where $\scf(G)$ can be arbitrary.  We begin with some results on the groups of linear characters of $\cH_{(\iota,\alpha,\beta)}$ and conclude with an antipode formula.

\subsection{Character groups} \label{CharacterGroupSection}

Recall, that a \emph{linear character} $\chi:\cH\rightarrow \CC$ is an algebra homomorphism, and the set of all such characters forms a group under the convolution product.  If $\Psi_\all:\ZZ_{\geq 0}\rightarrow \cH_{(\iota,\alpha,\beta)}$ is a function such that $\Psi_n\in \scf(G^{n-1})$, then we may define
\[
\begin{array}{rccl} \Psi_\all\rangle: & \cH_{(\iota,\alpha,\beta)} & \longrightarrow & \CC\\ & \gamma & \mapsto & \langle \gamma,\Psi_n\rangle \quad \text{if $\gamma\in \scf(G^{n-1})$.}
\end{array}\]
In fact, every linear character of $\cH_{(\iota,\alpha,\beta)}$ is of this form, although the converse does not hold.   The following lemma gives a characterization.
  
\begin{lemma} \label{AlgebraMorphismCondition}
Let $\Psi:\ZZ_{\geq 0} \rightarrow \cH_{(\iota,\alpha,\beta)}$ be a function such that $\Psi_n\in \scf(G^{n-1})$ and $\Psi_0=\chi^\emptyset$.  Then $\Psi_\bullet\rangle$ is an algebra homomorphism if and only if
for each $0\leq j\leq n$,
\[
\Def^{G^{n-1}}_{G^{j-1}\times G^{n-j-1}} (\iota_{(\{1,2,\ldots, j\},\{j+1,\ldots, n\})}^{(\{1,2,\ldots,n\})}\odot\Psi_n)=\Psi_j\otimes \Psi_{n-j}.
\]
\end{lemma}
\begin{proof}
On the one hand, if $\gamma\in \scf(G^{m-1})$ and $\chi\in \scf(G^{n-1})$, then
\begin{align*}
\langle \gamma\cdot \chi, \Psi_{m+n}\rangle&= \langle \iota_{(\{1,\ldots,m\},\{m+1,\ldots,m+n\})}^{(\{1,\ldots,m+n\})}\odot\Inf_{G^{m-1}\times G^{n-1}}^{G^{m+n-1}}(\gamma\otimes \chi),\Psi_{m+n}\rangle\\
&= \langle \gamma\otimes \chi,\Def_{G^{m-1}\times G^{n-1}}^{G^{m+n-1}}(\iota_{(\{1,\ldots,m\},\{m+1,\ldots,m+n\})}^{(\{1,\ldots,m+n\})}\odot\Psi_{m+n})\rangle.
\end{align*}
On the other hand,
\begin{equation*}
\langle \gamma, \Psi_{m}\rangle\langle \chi, \Psi_{n}\rangle = \langle \gamma\otimes \chi, \Psi_m\otimes \Psi_n\rangle.
\end{equation*}
Since $\gamma$ and $\chi$ were arbitrary, the result follows.
\end{proof}

From a representation theoretic point of view this construction is especially nice if each $\Psi_n$ is in fact the character of a module; in this case, if $\Psi_n$ is the trace of a $G^{n-1}$-module for all $n\geq 1$, then we say $\Psi_\all\rangle$ is \emph{supported by modules}.

One way to construct such a linear character is to fix $\psi\in \scf(G)$ and define  $\psi_\all\rangle$ by
\[
\psi_n=\underbrace{\psi\otimes \psi\otimes\cdots\otimes \psi}_{n-1\text{ terms}}\otimes \chi^{()}.
\]
Then
\begin{equation*}
\Def^{G^{n-1}}_{G^{j-1}\times G^{n-j-1}} (\iota_{(\{1,2,\ldots, j\},\{j+1,\ldots, n\})}^{(\{1,2,\ldots,n\})}\odot\psi_n)=\langle \psi,\iota\rangle\Big(\psi_j\otimes \psi_{n-j}\Big).
\end{equation*}
Thus, $\psi_\all\rangle$ is an algebra morphism if and only if $\langle\iota,\psi\rangle=1$.  In particular, $\psi\in \{\alpha,\beta\}$ always works for $\cH_{(\iota,\alpha,\beta)}$.

For $\Psi_\all\rangle$ an algebra homomorphism, $\tau\in \scf(G)$, and $\mu=(\mu_1,\ldots,\mu_\ell)\vDash n$, let
\[
(\Psi_\all)_\mu^\tau= \Psi_{\mu_1}\otimes\tau \otimes \Psi_{\mu_2}\otimes\tau \otimes \cdots \otimes \tau\otimes \Psi_{\mu_\ell},
\]
where we make use of the convention that
\[
\Psi_1\otimes \tau = \langle \Psi_1,\chi^{()}\rangle \tau
\]
as class functions of $\scf(G)$.  For $\rho\in \scf(G)$ and $\Gamma_\all\rangle$ another algebra homomorphism, let
\[
(\Psi_\all)\overset{\mu}{{}_\tau\shuffle_\rho}(\Gamma_\all)=(\Psi_{\mu_1}\otimes \tau)\otimes (\Gamma_{\mu_2}\otimes \rho)\otimes (\Psi_{\mu_3}\otimes \tau)\otimes (\Gamma_{\mu_4}\otimes \rho)\otimes \cdots\otimes \Omega_{\mu_\ell},\ \text{where}\ \Omega=\left\{\begin{array}{ll} \Psi & \text{if $\ell$ is odd},\\ \Gamma & \text{if $\ell$ is even.}\end{array}\right.
\]

\begin{theorem} \label{AlgebraMorphismGroup}
Let $\Psi_\all\rangle, \Gamma_\all\rangle$ be algebra morphisms. Then
\begin{enumerate}
\item[(a)]  $\Psi_\all\rangle\circ \Gamma_\all\rangle=(\Psi\circ\Gamma)_\all\rangle$ where
\[
(\Psi\circ\Gamma)_n=\sum_{\mu\vDash n} (\Psi_\all)\overset{\mu}{{}_\alpha\shuffle_\beta}(\Gamma_\all)+ (\Gamma_\all)\overset{\mu}{{}_\beta\shuffle_\alpha}(\Psi_\all).
\]
\item[(b)] $\Psi_\all\rangle^{-1}=(\Psi^{-1})_\all\rangle,$ where
\[
(\Psi^{-1})_n=\sum_{\mu\vDash n} (-\Psi_\all)_\mu^{\alpha+\beta}.
\]
\end{enumerate}
\end{theorem}
\begin{proof}
(a) We have that for $\psi\in \scf(G^{n-1})$, by Frobenius reciprocity,
\begin{align*}
(\Psi_\all\rangle\circ \Gamma_\all\rangle)(\psi)&=\sum_{A\subseteq \{1,2,\ldots,n\}} \langle\Down_{(A,\overline{A})}^{(\{1,2,\ldots,n\})}[\iota,\alpha,\beta](\psi),\Psi_{|A|}\times \Gamma_{|\overline{A}|} \rangle\\
&=\sum_{A\subseteq \{1,2,\ldots,n\} } \langle\psi,(\alpha\times \beta)_{(A,\overline{A})}^{(\{1,2,\ldots,n\})}\odot\Inf_{G^{n-1}_{\bc(A,\overline{A})}}^{G^{n-1}}\circ\Def_{G^{n-1}_{\bc(A,\overline{A})}}^{G^{n-1}_{\lc(A,\overline{A})}}(\iota_{(A,\overline{A})}\odot\Psi_{|A|}\times \Gamma_{|\overline{A}|}) \rangle,
\end{align*}
so
\[(\Psi\circ\Gamma)_n=\sum_{A\subseteq \{1,2,\ldots,n\} }(\alpha\times \beta)_{(A,\overline{A})}^{(\{1,2,\ldots,n\})}\odot\Inf_{G^{n-1}_{\bc(A,\overline{A})}}^{G^{n-1}}\circ\Def_{G^{n-1}_{\bc(A,\overline{A})}}^{G^{n-1}_{\lc(A,\overline{A})}}(\iota_{(A,\overline{A})}\odot\Psi_{|A|}\times \Gamma_{|\overline{A}|}). \]
As functions on 
\[
G^{n-1}_{\lc(A,\overline{A})}\cong G^{|A|-1}\times G^{|\overline{A}|-1}
\]
$\Psi_{|A|}$ evaluates on the first factor and $\Gamma_{|\overline{A}|}$ evaluates on the second.  When we deflate, $\iota_{(A,\overline{A})}$ multiplies each coordinate $i$ such that $i+1$ is in a different block of $(A,\overline{A})$ by $\iota$.  Due to this separation, we can use Lemma \ref{AlgebraMorphismCondition} to see that 
\[
\Def_{G^{n-1}_{\bc(A,\overline{A})}}^{G^{n-1}_{\lc(A,\overline{A})}}(\iota_{(A,\overline{A})}\odot\Psi_{|A|}\times \Gamma_{|\overline{A}|})=\Psi_{\mu_A}\otimes \Gamma_{\mu_{\overline{A}}},
\]
where $\mu_A$ is the integer composition whose parts lengths record the length of the sequences with consecutive integers in $A$ and $\mu_{\overline{A}}$  is the integer composition whose parts lengths record the length of the sequences with consecutive integers in $\overline{A}$ .  For example, if $A=\{1,4,5,7,8,9\}\subseteq \{1,2,\ldots, 15\}$, then
\[
\mu_A=
\begin{tikzpicture}[scale=.3,baseline=.3cm]
\foreach \x/\y/\z in {0/2/1,0/1/4,1/1/5,1/0/7,2/0/8,3/0/9}
	{\draw[red] (\x,\y) +(-.5,-.5) rectangle ++(.5,.5);
	\node at (\x,\y) {$\scscs{\z}$};}
\end{tikzpicture}
\quad
\text{and}
\quad
\mu_{\overline{A}}=
\begin{tikzpicture}[scale=.3,baseline=.3cm]
\foreach \x/\y/\z in {0/2/2,1/2/3,1/1/6,1/0/10,2/0/11,3/0/12,4/0/13,5/0/14,6/0/15}
	{\draw[blue] (\x,\y) +(-.5,-.5) rectangle ++(.5,.5);
	\node at (\x,\y) {$\scscs{\z}$};}
\end{tikzpicture}\ .
\]
However, as we inflate to $G^{n-1}$ the relative order of the parts is preserved, so something like
\[
\begin{tikzpicture}[scale=.3,baseline=.6cm]
\foreach \x/\y/\z in {0/5/1,1/3/4,2/3/5,2/1/7,3/1/8,4/1/9}
	{\draw[red] (\x,\y) +(-.5,-.5) rectangle ++(.5,.5);
	\node at (\x,\y) {$\scscs{\z}$};}
\foreach \x/\y/\z in {0/4/2,1/4/3,2/2/6,4/0/10,5/0/11,6/0/12,7/0/13,8/0/14,9/0/15}
	{\draw[blue] (\x,\y) +(-.5,-.5) rectangle ++(.5,.5);
	\node at (\x,\y) {$\scscs{\z}$};}
\foreach \x/\y/\l in {10/0/{\Gamma_6},5/1/{\Psi_3},3/2/{\Gamma_1},3/3/{\Psi_2},2/4/{\Gamma_2},1/5/{\Psi_1}}
	\node at (\x+.25,\y) {$\scscs\l$};
\end{tikzpicture}\ .
\]
Now the $(\alpha\times \beta)_{(A,\overline{A})}^{(\{1,2,\ldots,n\})}$ term adds an $\alpha$ whenever we transition from $A$ to $\overline{A}$ and a $\beta$ whenever we go from $\overline{A}$ to $A$. Result (a) follows as we sum over all $A$.

(b) Let $\Gamma_\all$ be given by 
\[\Gamma_n=\sum_{\mu\vDash n} (-\Psi_\all)_\mu^{\alpha+\beta}.\]
We want to show that $\Gamma_\all\rangle=\Psi_\all\rangle^{-1}$ or $(\Psi_\all\circ\Gamma_\all)_n=0$ for all $n\geq 1$.  By (a),
\begin{equation*}
(\Psi_\all\circ \Gamma_\all)_n=\sum_{\mu\vDash n} (\Psi_\all)\overset{\mu}{{}_\alpha\shuffle_\beta}(\Gamma_\all)+ (\Gamma_\all)\overset{\mu}{{}_\beta\shuffle_\alpha}(\Psi_\all)
\end{equation*}
If we expand the $\Gamma$-terms according to its definition, all the summands are of the form.
\begin{equation}\label{ProductSummands}
\pm\Psi_{\nu_1}\otimes x_1\otimes \Psi_{\nu_2}\otimes x_2\otimes \cdots \otimes x_{\ell-1} \otimes\Psi_{\nu_\ell}
\end{equation}
where $\nu=(\nu_1,\ldots,\nu_\ell)\vDash n$ and $x_1,\ldots,x_{\ell-1}\in \{\alpha,\beta,\alpha+\beta\}$.  However, for a fixed composition $\nu$ not all sequences in $\{\alpha,\beta,\alpha+\beta\}^{\ell-1}$ can occur.  Since $\Psi_k$ and $\Gamma_l$ terms alternate, any $\alpha$ in the sequence must be followed by either an $\alpha+\beta$ or a $\beta$.  Similarly, any $\beta$ must be followed by an $\alpha$. The $\alpha+\beta$ could be followed by either a $\beta$ or an $\alpha+\beta$.  Define
\[
S_{\ell}=\left\{(x_1,\ldots,x_{\ell-1})\in \{\alpha,\beta,\alpha+\beta\}^{\ell-1}\ \bigg|\ \begin{array}{@{}l@{}} x_i\in\{\alpha,\alpha+\beta\}\text{ implies } x_{i+1}\in \{\beta,\alpha+\beta\},\\  x_i=\beta \text{ implies }  x_{i+1}=\alpha\end{array}\right\}
\]
This set is in bijection with the summands of the form (\ref{ProductSummands}).  For example,
\begin{align*}
(\alpha,\beta,\alpha,\alpha+\beta&,\alpha+\beta,\beta,\alpha,\alpha+\beta)\\
&\updownarrow\\
\overbrace{\Psi_{\nu_1}}^{{\color{red} \Psi_{\nu_1}}}\otimes \alpha \otimes \overbrace{ (-\Psi_{\nu_2})}^{{\color{blue} \Gamma_{\nu_2}}}\otimes \beta\otimes \overbrace{\Psi_{\nu_3}}^{{\color{red} \Psi_{\nu_3}}}\otimes \alpha
\otimes &\overbrace{(-\Psi_{\nu_4})\otimes (\alpha+\beta)\otimes (-\Psi_{\nu_5})\otimes (\alpha+\beta)\otimes (-\Psi_{\nu_6})}^{{\color{blue}\Gamma_{\nu_4+\nu_5+\nu_6}}}\\
\otimes & \beta  \otimes  \underbrace{\Psi_{\nu_7}}_{\color{red}\Psi_{\nu_7}}\otimes \alpha\otimes \underbrace{(-\Psi_{\nu_8})\otimes (\alpha+\beta)\otimes (-\Psi_{\nu_{9}})}_{\color{blue}\Gamma_{\nu_8+\nu_9}}
\end{align*}
The sign of the summand is also determined by
\[
\sign(x_1,\ldots,x_{\ell-1})=(-1)^{\ell-1-\#\{i\mid x_i=\alpha\}-\#\{x_{\ell-1}=\beta\}}.
\]
Consider the sign-reversing bijections 
\[
\begin{array}{ccc} \left\{(x_1,\ldots,x_{\ell-1})\in S_\ell\ \bigg|\ \begin{array}{@{}l@{}} x_j=\beta,x_{j+1}=\alpha,\\
x_k=\alpha+\beta,k>j+1\end{array} \right\} & \overset{\toggle_j}{\longleftrightarrow} &  \{(x_1,\ldots,x_{\ell-1})\in S_\ell\mid x_k=\alpha+\beta,k\geq j\}\\
(x_1,\ldots,x_{j-1},\beta,\alpha,x_{j+1},\ldots,x_{\ell-1}) & \mapsto &  (x_1,\ldots,x_{j-1},\alpha+\beta,\alpha+\beta,x_{j+1},\ldots,x_{\ell-1}) \end{array}
\]
in $S_\ell$. Note these functions are well-defined since both $\beta$ and $\alpha+\beta$ may only be preceded by $\alpha$ or $\alpha+\beta$. Now by $\toggle_{\ell-1}$,
\begin{align*}
\sum_{(x_1,\ldots,x_{\ell-1})\in S_\ell\atop x_{\ell-1}=\beta} &\sign(x_1,\ldots, x_{\ell-1})x_1\otimes \cdots\otimes x_{\ell-1}+\sum_{(x_1,\ldots,x_{\ell-1})\in S_\ell\atop x_{\ell-1}=\alpha+\beta} \sign(x_1,\ldots, x_{\ell-1})x_1\otimes \cdots\otimes x_{\ell-1}\\
&=\sum_{(x_1,\ldots,x_{\ell-1})\in S_\ell\atop x_{\ell-1}=\alpha+\beta} \sign(x_1,\ldots, x_{\ell-1})x_1\otimes \cdots \otimes x_{\ell-2}\otimes \alpha.
\end{align*}
Thus,
\begin{align*}
&\sum_{(x_1,\ldots,x_{\ell-1})\in S_\ell} \sign(x_1,\ldots, x_{\ell-1})x_1\otimes \cdots\otimes x_{\ell-1}\\
&=\sum_{(x_1,\ldots,x_{\ell-1})\in S_\ell\atop x_{\ell-2}=\beta, x_{\ell-1}=\alpha} \sign(x_1,\ldots, x_{\ell-1})x_1\otimes \cdots\otimes x_{\ell-1}+\sum_{(x_1,\ldots,x_{\ell-1})\in S_\ell\atop x_{\ell-1}=\alpha+\beta} \sign(x_1,\ldots, x_{\ell-1})x_1\otimes \cdots \otimes x_{\ell-2}\otimes \alpha\\
&=\hspace{-.5cm} \sum_{(x_1,\ldots,x_{\ell-1})\in S_\ell\atop x_{\ell-3}=\beta,x_{\ell-2}=\alpha, x_{\ell-1}=\alpha+\beta} \hspace{-1cm}\sign(x_1,\ldots, x_{\ell-1}) x_1\otimes \cdots\otimes x_{\ell-2}\otimes \alpha+\hspace{-.75cm} \sum_{(x_1,\ldots,x_{\ell-1})\in S_\ell\atop x_{\ell-2}=x_{\ell-1}=\alpha+\beta} \hspace{-.75cm} \sign(x_1,\ldots, x_{\ell-1})x_1\otimes \cdots \otimes x_{\ell-3}\otimes \alpha\otimes \alpha,
\end{align*}
by applying $\toggle_{\ell-2}$.  Continue to iterate with $\toggle_{\ell-3}$, $\toggle_{\ell-4}$,..., $\toggle_1$,  until we get to
\begin{align*}
&\sum_{(x_1,\ldots,x_{\ell-1})\in S_\ell} \sign(x_1,\ldots, x_{\ell-1})x_1\otimes \cdots\otimes x_{\ell-1}\\
&=\hspace{-.5cm}\sum_{(x_1,\ldots,x_{\ell-1})\in S_\ell\atop x_{1}=\alpha, x_2=x_3=\cdots=x_{\ell-1}=\alpha+\beta} \hspace{-.75cm} \sign(x_1,\ldots, x_{\ell-1}) \alpha\otimes \cdots \otimes \alpha+
\hspace{-.5cm}\sum_{(x_1,\ldots,x_{\ell-1})\in S_\ell\atop x_{1}=x_{2}=\cdots =x_{\ell-1}=\alpha+\beta} \hspace{-.5cm}\sign(x_1,\ldots, x_{\ell-1})\alpha\otimes \cdots \otimes \alpha\\
&=0,
\end{align*}
since these last two summands have opposite sign.  This implies that all the summands of the form (\ref{ProductSummands}) sum to zero and so $(\Psi_\all\circ\Gamma_\all)_n=0$ for $n\geq 1$.
\end{proof}

\begin{corollary}  Suppose $\alpha$ and $\beta$ are class functions that are traces of modules.  
If $\Psi_\all\rangle$ and $\Gamma_\all\rangle$ are algebra homomorphisms of $\cH_{(\iota,\alpha,\beta)}$ supported by modules, then $(\Psi\circ\Gamma)_\all\rangle$ is an algebra homomorphism supported by modules.
\end{corollary}

\begin{corollary}
If $\alpha=\beta$, then $\alpha_\all\rangle$ is an odd linear character of $\cH_{(\iota,\alpha,\beta)}$.
\end{corollary}
\begin{proof} Note that since $\alpha=\beta$, for $\mu\vDash n$,
\[
(\alpha_\bullet)_\mu^{\alpha+\beta}=(\alpha_\bullet)_\mu^{2\alpha}=2^{\ell(\mu)-1} (\alpha_\bullet)_{(n)}^\iota.
\]
Thus, by Theorem \ref{AlgebraMorphismGroup} (b),
\begin{align*}
(\alpha_\all^{-1})_n
&=\sum_{\mu\vDash n} (-1)^{\ell(\mu)}2^{\ell(\mu)-1} (\alpha_\all)^\iota_{(n)}\\
&=\Big(\sum_{A\subseteq \{1,2,\ldots,n-1\}} (-1)^{|A|+1}2^{|A|} \Big)(\alpha_\all)^\iota_{(n)}\\
&=\Big(-\sum_{A\subseteq \{1,2,\ldots,n-1\}} (-2)^{|A|} \Big)(\alpha_\all)^\iota_{(n)}\\
&=-(-2+1)^{n-1}(\alpha_\all)^\iota_{(n)}\\
&=(-1)^n(\alpha_\all)^\iota_{(n)},
\end{align*}
as desired.
\end{proof}
\subsection{Antipode} \label{Antipode}

For $\gamma_1,\ldots, \gamma_{n-1}\in \scf(G)$ and $\gamma_n=\chi^{()}$,
\[
S(\gamma_1\otimes \cdots \otimes \gamma_{n}) = \sum_{\underline{A}\vDash \{1,2,\ldots,n\}} (-1)^{\ell(\underline{A})-1}\Inf^{(\{1,2,\ldots,n\})}_{w_{\underline{A}}(\underline{A})}[\iota]\circ w_{\underline{A}} \circ\Down^{(\{1,2,\ldots,n\})}_{\underline{A}}[\iota,\alpha,\beta](\gamma_1\otimes \cdots \otimes \gamma_{n}),
\]
where $w_{\underline{A}}$ is the straightening permutation with minimal crossings that sends 
\[A_1\mapsto \{1,2,\ldots, |A_1|\},A_2\mapsto \{|A_1|+1,\ldots, |A_1|+|A_2|\},\ldots, A_\ell\mapsto \{n-|A_\ell|+1,\ldots, n\}.\]
Specifically, for $j\in A_k$,
\[w_{\underline{A}}(j)=\#\{i\leq j\mid i\in A_k\} +\sum_{i<k} |A_i|.
\]
Note that 
\begin{align*}
\Inf^{(\{1,2,\ldots,n\})}_{w_{\underline{A}}(\underline{A})}[\iota]\circ w_{\underline{A}}\circ&\Down^{(\{1,2,\ldots,n\})}_{\underline{A}}[\iota,\alpha,\beta](\gamma_1\otimes \cdots \otimes \gamma_{n})_{w_{\underline{A}}(j)}\\
&=\left\{\begin{array}{ll}
\gamma_j & \text{if $\bc(\underline{A})_j=1$,}\\
\langle\gamma_j,\alpha\rangle\iota & \text{if $\bc(\underline{A})_j=0$, $\llc(\underline{A})_j=1$,}\\ 
\langle\gamma_j,\beta\rangle\iota & \text{if $\bc(\underline{A})_j=0$, $\llc(\underline{A})_j=0$, $w_{\underline{A}}(j)\neq n$},\\
\langle\gamma_j,\beta\rangle\chi^{()} & \text{if $\bc(\underline{A})_j=0$, $\llc(\underline{A})_j=0$, $w_{\underline{A}}(j)= n$},\\
\iota & \text{if $j=n\neq w_{\underline{A}}(j)$,}\\
\chi^{()}& \text{if $j=n= w_{\underline{A}}(j)$.}\\
\end{array}\right.
\end{align*}
Consider the case where we have $j<k$ and 
\[ w_{\underline{A}}(k)=w_{\underline{A}}(j)+1,\]
or
\[
\begin{tikzpicture}[baseline=.8cm]
\node at (2,1.5) {$\otimes\gamma_j\otimes \gamma_{j+1}\otimes \cdots \otimes \gamma_{k}$};
\node at (0,0) {$\otimes \cdot_{w_{\underline{A}}(j)} \otimes \cdot_{w_{\underline{A}}(k)}\otimes \cdots $};
\draw (0.25,1.35) -- (-1.45,.25);
\draw (3,1.35) -- (-.15,.25);
\end{tikzpicture}\ .
\]
Then we could either have $j$ and $k$ in the same block of $\underline{A}$ or not.   If $k\neq j+1$, then we get the same answer for the $w_{\underline{A}}(j)$ coordinate whether $j$ and $k$ are in the same block of $\underline{A}$ or not.   However, these two choices do have a different sign in the sum.  If $k=j+1$, we get the same answer only when $\gamma_j=\iota$.  With this in mind, a \emph{toggle point} in $\underline{A}\vDash \{1,2,\ldots, n\}$ is an element $1\leq j\leq n$ such that if $j\in A_i$, then either
\begin{description}
\item[split toggle.] $j$ is maximal in $A_i$ and $\min\{A_{i+1}\}>j+1$, or
\item[fused toggle.] $j$ is not maximal in $A_i$ and $j+1\notin A_i$.
\end{description}
Define
\begin{align*}
\TF_n&=\{\underline{A}\vDash \{1,2,\ldots, n\}\mid \underline{A}\text{ has no toggle point}\}
\end{align*}
\begin{lemma}\hfill

\begin{enumerate}
\item[(a)] $|\TF_n|=3^{n-1}$,
\item[(b)] \[
S(\gamma_1\otimes \cdots \otimes \gamma_{n}) = \sum_{\underline{A}\in \TF_n} (-1)^{\ell(\underline{A})-1}\Inf^{(\{1,2,\ldots,n\})}_{w_{\underline{A}}(\underline{A})}[\iota]\circ w_{\underline{A}} \circ\Down^{(\{1,2,\ldots,n\})}_{\underline{A}}[\iota,\alpha,\beta](\gamma_1\otimes \cdots \otimes \gamma_{n}).
\]
\end{enumerate}
\end{lemma}
\begin{proof}
(a) To construct a generic element of $\TF_n$, make the following choices:
\begin{description}
\item[Step 1.]  Place 2 in the same block as 1,  start  the next block with 2, or start a block in the beginning with 2.
\item[Step 2.] Place 3 in the same block as 2, start the next block with 3, or start a block in the beginning with 3.
\item[$\vdots$]
\item[Step $n-1$.] Place $n$ in the same block as $n-1$, start the next block with $n$, or start a block in the beginning with $n$.
\end{description}
Constructing the first few examples is as follows:
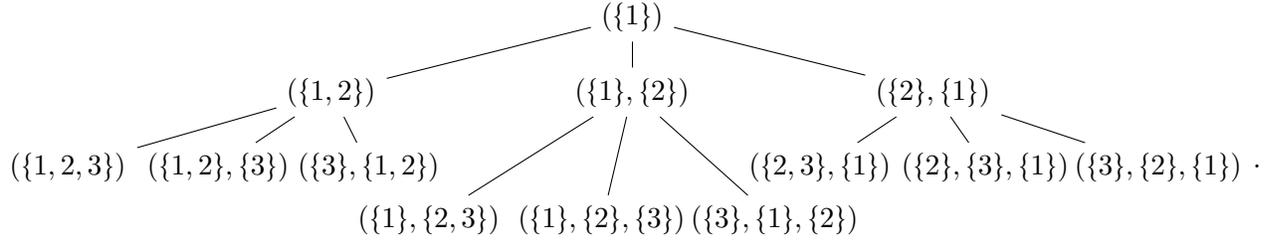
\[\begin{tikzpicture}[baseline=3cm]
\node (1) at (0,5) {$(\{1\})$};
\node (11) at (-4,4) {$(\{1,2\})$};
\node (111) at (-7.5,3) {$(\{1,2,3\})$};
\node (112) at (-5.5,3) {$(\{1,2\},\{3\})$};
\node (113) at (-3.5,3) {$(\{3\},\{1,2\})$};
\node (12) at (0,4) {$(\{1\},\{2\})$};
\node (121) at (-2.7,2.3) {$(\{1\},\{2,3\})$};
\node (122)  at (-.4,2.3) {$(\{1\},\{2\},\{3\})$};
\node (123) at (1.9,2.3)  {$(\{3\},\{1\},\{2\})$};
\node (13) at (4,4) {$(\{2\},\{1\})$};
\node (131) at (2.5,3) {$(\{2,3\},\{1\})$};
\node (132) at (4.7,3) {$(\{2\},\{3\},\{1\})$};
\node (133) at (7,3) {$(\{3\},\{2\},\{1\})$};
\foreach \x in {1,2,3}
	{\draw (1) -- (1\x);
	\draw (11) -- (11\x);
	\draw (12) -- (12\x);
	\draw (13) -- (13\x);}
\end{tikzpicture} .\]
In each case, we have 3 choices, giving the total.

(b) We begin by sorting the set compositions of $\{1,2,\ldots,n\}$ that have toggles into two subsets:
\begin{align*}
\Splt & =\{\underline{A}\vDash \{1,2,\ldots,n\} \mid \underline{A}\notin\TF_n,\text{ the minimal toggle is split}\}\\
\Fusd&= \{\underline{A}\vDash \{1,2,\ldots,n\} \mid \underline{A}\notin\TF_n,\text{ the minimal toggle is fused}\}.
\end{align*}
By changing the minimal toggle from split to fused (fusing the corresponding blocks), we obtain a sign reversing bijection
\[
\Splt\longleftrightarrow 
\Fusd,
\]
where, for example,
\[(\{5,6,7\},\{1,2,\underline{3\},\{8},9\},\{4,10\})\leftrightarrow (\{5,6,7\},\{1,2,\underline{3 , 8},9\},\{4,10\}).\]
However, because $\llc$, $\bc$ and $\lc$ are the same between the two cases, the terms cancel in the sum.
\end{proof}

For $\mu\vDash n$, let
\[
\overline{\mu_{i}}=\mu_1+\mu_2+\cdots+\mu_i,
\]
and let $\overline{\mu_0}=0$.
\begin{theorem}
\[
S(\gamma_1\otimes \cdots \otimes \gamma_{n})=\sum_{\mu\vDash n}\Big(\prod_{j=1}^{\ell-1} \langle \gamma_{\overline{\mu_j}},\beta\rangle\Big) (-\gamma_{\mu}^{(\ell)})\otimes \iota\otimes  (-\gamma_{\mu}^{(\ell-1)})\otimes \iota\otimes \cdots\otimes \iota\otimes (-\gamma_{\mu}^{(1)}),
\]
where
\[
\gamma_{\mu}^{(j)}=(\gamma_{\overline{\mu_{j-1}}+1}-\langle\gamma_{\overline{\mu_{j-1}}+1},\alpha\rangle \iota)\otimes  \cdots \otimes (\gamma_{\overline{\mu_{j}}-1}-\langle \gamma_{\overline{\mu_{j}}-1},\alpha\rangle\iota)\in \scf(G^{\mu_j-1}).
\]
\end{theorem}
\begin{remark}
This formula is generically cancellation free. That is, if $\gamma_1,\ldots,\gamma_n$ are linearly independent and do not include a multiple of $\iota$, then there is no cancellation.
\end{remark}

\begin{proof}
We begin an equivalence relation on $\TF_n$.  Let $\underline{A}\sim \underline{B}$ if $w_{\underline{A}}^{-1}=w_{\underline{B}}^{-1}$.  This means that if one reads off the integers in increasing order in the order of the parts, then the resulting words are the same ($w_{\underline{A}}^{-1}$ in 1-line notation), so for example
\begin{equation}\label{PermutationExample}
\begin{split} \underline{A}&=(\{7,8\},\{9,10\},\{3,4,5,6\}, \{1\},\{2\})\\ \underline{B}&=(\{7\},\{8, 9,10\},\{3,4\},\{5,6\}, \{1,2\})
\end{split}
\quad \text{have}\quad w_{\underline{A}}^{-1}=(7,8,9,10,3,4,5,6,1,2)=w_{\underline{B}}^{-1}.
\end{equation}
Note that 
\begin{align*}\sum_{w_{\underline{A}}^{-1}=(1,2,\ldots,n)}& (-1)^{\ell(\underline{A})}\Inf^{(\{1,2,\ldots,n\})}_{w_{\underline{A}}(\underline{A})}[\iota]\circ w_{\underline{A}} \circ\Down^{(\{1,2,\ldots,n\})}_{\underline{A}}[\iota,\alpha,\beta](\gamma_1\otimes \cdots \otimes \gamma_{n})\\
&= \sum_{w_{\underline{A}}^{-1}=(1,2,\ldots,n)}   (-1)^{\ell(\underline{A})}  \gamma_{A_1}\otimes \langle \gamma_{|A_1|},\alpha\rangle \iota \otimes  \gamma_{A_2}\otimes \langle \gamma_{|A_1\cup A_2|},\alpha\rangle \iota \otimes \cdots \otimes  \gamma_{A_{\ell}},
\end{align*}
where
\[
\gamma_{A_j}=\gamma_{i_1}\otimes \gamma_{i_2}\otimes\cdots\otimes \gamma_{i_{k-1}}\quad \text{if}\quad A_j=\{i_1<i_2<\cdots<i_k\}.
\]
By comparing the coefficients of the monomial terms we see that
\begin{align*}
 \sum_{w_{\underline{A}}^{-1}=(1,2,\ldots,n)}  & (-1)^{\ell(\underline{A})}  \gamma_{A_1}\otimes \langle \gamma_{|A_1|},\alpha\rangle \iota \otimes  \gamma_{A_2}\otimes \langle \gamma_{|A_1\cup A_2|},\alpha\rangle \iota \otimes \cdots \otimes  \gamma_{A_{\ell}}\\
 &=-(\gamma_1- \langle \gamma_1,\alpha\rangle \iota)\otimes(\gamma_2- \langle \gamma_2,\alpha\rangle \iota) \otimes \cdots \otimes (\gamma_{n-1}- \langle \gamma_{n-1},\alpha\rangle \iota)\otimes \chi^{()}\\
 &=-\gamma_{(n)}^{(1)}.
\end{align*}
Now fix $w_{\underline{A}}^{-1}$ with $\underline{A}\in \TF_n$, and let $\mu$ be the composition that gives the lengths of the consecutive rising sequences in $\underline{A}$; in (\ref{PermutationExample}), $\mu=(4,4,2)$.  Then
\begin{align*}
& \sum_{w_{\underline{B}}^{-1}=w_{\underline{A}}^{-1}}   (-1)^{\ell(\underline{B})}\Inf^{(\{1,2,\ldots,n\})}_{w_{\underline{B}}(\underline{B})}[\iota]\circ w_{\underline{B}} \circ\Down^{(\{1,2,\ldots,n\})}_{\underline{B}}[\iota,\alpha,\beta](\gamma_1\otimes \cdots \otimes \gamma_{n})\\
 &=\sum_{\underline{B}\vDash\{\mu_1+1,\ldots, n\}\atop w_{\underline{B}}^{-1}=w_{(A_1,\ldots,A_{\ell-1})}^{-1}} (-1)^{\ell(\underline{B})}\Inf^{(\{1,2,\ldots,n\})}_{w_{\underline{B}}(\underline{B})}[\iota]\circ w_{\underline{B}} \circ\Down^{(\{1,2,\ldots,n\})}_{\underline{B}}[\iota,\alpha,\beta](\gamma_{\mu_1+1}\otimes \cdots \otimes \gamma_{n})\\
 &\hspace*{1cm} \otimes \hspace{-.3cm} \sum_{\underline{C}\vDash\{1,2,\ldots, \mu_1\}} \hspace{-.3cm} (-1)^{\ell(\underline{C})}\Inf^{(\{1,2,\ldots,n\})}_{w_{\underline{C}}(\underline{C})}[\iota]\circ w_{\underline{C}} \circ\Down^{(\{1,2,\ldots,n\})}_{\underline{C}}[\iota,\alpha,\beta](\gamma_{1}\otimes \cdots\otimes\gamma_{\mu_1-1} \otimes \langle\gamma_{\mu_1},\beta\rangle\chi^{()}) \\
  &=\hspace{-.6cm}\sum_{\underline{B}\vDash\{\mu_1+1,\ldots, n\}\atop w_{\underline{B}}^{-1}=w_{(A_1,\ldots,A_{\ell-1})}^{-1}} \hspace{-.7cm}(-1)^{\ell(\underline{B})}\langle\gamma_{\mu_1},\beta\rangle\Inf^{(\{1,2,\ldots,n\})}_{w_{\underline{B}}(\underline{B})}[\iota]\circ w_{\underline{B}} \circ\Down^{(\{1,2,\ldots,n\})}_{\underline{B}}[\iota,\alpha,\beta](\gamma_{\mu_1+1}\otimes \cdots \otimes \gamma_{n}) \otimes(-\gamma_{\mu}^{(1)}). 
\end{align*}
Iterate to get the result.
\end{proof}

If we return to the situation of Lemma \ref{ExplicitNSymIsos}, we have the following consequences.

\begin{corollary} \label{GenericAntipode}
Let $\scf(G)=\CC\spanning\{\mathds{1},\reg\}$.  Then
\begin{description}
\item[Case $\alpha=\beta$.]  If $\tau\in \scf(G)$ is a nonzero element such that $\langle\tau,\beta\rangle=0$ and $\CC\spanning\{\tau,\iota\}=\scf(G)$, then
\[
S((\tau_\all)_{(n)}^\iota)=-(\tau_\all)_{(n)}^\iota.
\]
\item[Case $\alpha\neq\beta$.]  If $\tau\in \scf(G)$ satisfies $\langle \tau,\alpha\rangle=1$ and $\langle \tau,\beta\rangle=0$, then
\[
S((\tau_\all)_{(n)}^\iota)=-((\tau-\iota)_\all)_{(n)}^\iota.
\]
\end{description}
\end{corollary}

Applying to $\NSym$ we do not get new results, but they are easy to prove.

\begin{corollary}
Let $\scf(G)=\CC\spanning\{\mathds{1},\reg\}$.  Then
\begin{description}
\item[Case $\alpha=\beta=\iota=\mathds{1}$.]  For $\mu\vDash n$, 
\[
S(\chi^{\binaryzero(\mu)})=(-1)^{\ell(\mu)} \chi^{\binaryzero(\mu_\ell)} \chi^{\binaryzero(\mu_{\ell-1})}\cdots \chi^{\binaryzero(\mu_{1})}.
\]
\item[Case $\iota=\reg,\alpha=\mathds{1},\beta=(\reg-\mathds{1})^*$.] 
 For $\mu\vDash n$, 
\[
S\Big(\Ind_{G_{\binaryone(\mu)}^{n-1}}^{G^{n-1}}(\mathds{1})\Big)=\sum_{\nu\text{ refines }  (\mu_\ell,\ldots,\mu_2,\mu_1)} (-1)^{\ell(\nu)}\Ind_{G_{\binaryone(\nu)}^{n-1}}^{G^{n-1}}(\mathds{1}).
\]
\end{description}
\end{corollary}

\begin{proof}
(a) Note that $\chi^{\binaryzero(\mu)}=(\reg-\mathds{1})_\mu^{\mathds{1}}$, so this follows from Corollary \ref{GenericAntipode} with $\tau=\reg-\mathds{1}$.

(b) Note that $\Ind_{G_{\binaryone(\mu)}^{n-1}}^{G^{n-1}}(\mathds{1})=(\mathds{1})_\mu^{\reg}$, so we will use Corollary \ref{GenericAntipode} with $\tau=\mathds{1}$.  In the case $\mu=(n)$, we have
\begin{align*}
S\Big((\mathds{1})_{(n)}^{\reg}\Big) & = -(\mathds{1}-\reg)_{(n)}^{\reg}\\
&= -\sum_{\nu\vDash n} (-1)^{\ell(\nu)-1} (\mathds{1})^{\reg}_\nu\\
&=\sum_{\nu\vDash n} (-1)^{\ell(\nu)} (\mathds{1})^{\reg}_\nu.
\end{align*}
The result now follows from the multiplicativity of the antipode.
\end{proof}

\begin{remark}
One can also prove the formula for the ribbon basis, but it does not fall directly out of the formula in the same way (there is a lot more cancellation involved).
\end{remark}

\end{document}